\documentclass[preprint,12pt]{elsarticle}
\usepackage[a4paper, total={6.5in, 8in}]{geometry}

\usepackage{amsmath}
\usepackage{graphicx} 
\usepackage{amssymb}
\usepackage{subcaption}
\newtheorem{remark}{Remark}

\newcommand{\dt}[1]{\frac{\text{d}#1}{\text{d}t}}
\newcommand{\dd}[2]{\frac{\text{d}#1}{\text{d}#2}} 
\newcommand{\ddq}[3]{\frac{\text{d}^{#1}#2}{\text{d}#3^{#1}}} 

\usepackage{xcolor}

 \journal{Journal}

\begin{document}

\begin{frontmatter}

\title{Deep smoothness WENO scheme for two-dimensional hyperbolic conservation laws: A deep learning approach for learning smoothness indicators}

\author{Tatiana Kossaczká$^{{\dagger},*}$, Ameya D. Jagtap$^{\ddagger}$, Matthias Ehrhardt$^{\dagger}$}
\cortext[mycorrespondingauthor]{Corresponding author emails: kossaczka@uni-wuppertal.de (T. Kossaczká), ameya\_jagtap@brown.edu (A. D. Jagtap), ehrhardt@uni-wuppertal.de (M. Ehrhardt)}

\address{$^{\dagger}$~Chair of Applied and Computational Mathematics, Bergische Universität Wuppertal, Gaußstrasse 20, Wuppertal, 42119, Germany}
\address{$^{\ddagger}$~Division of Applied Mathematics, Brown University, 182 George Street, Providence, RI, 02912, USA}

\begin{abstract}
In this paper, we introduce an improved version of the fifth-order weighted essentially non-oscillatory (WENO) shock-capturing scheme by incorporating deep learning techniques. 
The established WENO algorithm is improved 
by training a compact neural network to adjust the smoothness indicators within the WENO scheme. 
This modification enhances the accuracy of the numerical results, particularly near abrupt shocks. 
Unlike previous deep learning-based methods, no additional post-processing steps are necessary for maintaining consistency. 
We demonstrate the superiority of our new approach using several examples from the literature for the two-dimensional Euler equations of gas dynamics. 
Through intensive study of these test problems, which involve various shocks and rarefaction waves, the new technique is shown to outperform traditional fifth-order WENO schemes, especially in cases where the numerical solutions exhibit excessive diffusion or overshoot around shocks. 
\end{abstract}




\end{frontmatter}


\section{Introduction}\label{sec:S1}
It has long been a challenge to adequately simulate complex flow problems using numerical methods. 
Recently, this has been further improved using machine learning techniques. 
As an example, in \cite{RAISSI2019686,jagtap2022physics,jagtap2022deep}, the concept of physics-informed neural networks (PINNs) for the solution of complex fluid flow problems was proposed, which seamlessly combines the data and the mathematical models; see  \cite{RAISSI2019686, jagtap2020conservative, jagtap2020extended, shukla2021parallel, PENWARDEN2023112464,de2022error} for more details.
Similarly, a new method using a U-Net-like convolutional neural network (CNN) along with established finite difference discretization techniques was proposed to learn approximate solutions for the NSE without the need for parameterization \cite{grimm2023learning}. 
Also, recently, a framework called \textit{local transfer function analysis} (LTA) for optimizing numerical methods for convection problems using a graph neural network (GNN) was proposed \cite{drozda2023learning}.

The work \cite{mao2020physics} investigated the use of PINNs to approximate the hyperbolic Euler equations of gas dynamics. 
The Euler equations and initial and boundary conditions are used to create a loss function that solves scenarios with smooth solutions and those with discontinuities. 
%
Next, in \cite{jagtap2020conservative}, a novel approach, called \textit{conservative PINNs}, for solving nonlinear conservation laws, such as the compressible Euler equations, was presented.  
%
%
%
In the recent paper \cite{van2023accelerating}, another novel approach has been proposed where machine learning improves finite-difference-based approximations of PDEs while maintaining high-order convergence through node refinement.

This research area is also the context of our work. 
Recently, improvements to the standard finite difference methods (FDMs) have been developed \cite{kossaczka2023}. By adding a small convolutional neural network, the solutions of the standard PDEs are improved, while the convergence and consistency properties of the original methods are preserved.
We aim to further improve modern FDMs, such as WENO schemes, for nonlinear hyperbolic systems using machine learning. 
For this type of PDEs, it is known that discontinuities (shocks) can occur despite initial smoothness, which makes specialized numerical methods mandatory.
Therefore, the focus of our attention is on the behavior of numerical solutions in the vicinity of shocks.

To better frame our current work, let us very briefly sketch the historical development of WENO schemes.
Crandall and Majda \cite{Crandall} introduced \textit{monotone schemes} in 1980 that maintain stability and satisfy entropy conditions, but are only exactly first order due to Godunov's theorem. 
Next, \textit{shock-capturing schemes} were developed to accurately handle shocks and gradients without excessive diffusion \cite{Harten}. 
The \textit{essentially non-oscillatory (ENO) schemes} \cite{Harten87} were outstanding, achieving high accuracy in smooth regions and effective shock resolution using smoothness indicators, e.g.\ \cite{jiang1996, shu1998}. 
Extensions such as the \textit{Hermite WENO (HWENO)} schemes \cite{qiu2004hermite, qiu2005hermite} and \textit{hybrid methods} \cite{pirozzoli2002conservative, hill2004hybrid} were introduced for higher accuracy and efficiency. 
A gas-kinetic theory-based KWENO scheme was proposed in \cite{jagtap2020kinetic} for hyperbolic conservation laws.
Moreover, further modifications of WENO scheme have been developed, e.g. \cite{Henrick2005, Borges2008, Castro2011, Ha2013, zhu2016new, rathan2020l1}.

Neural networks approximated the solutions of PDEs and improved numerical methods for PDEs. 
While the data-driven approach is promising for improving modern numerical methods, it is always important to maintain a balance between new data-driven insights and established mathematical structures, i.e., the basic numerical scheme (here based on physical principles), e.g., for hyperbolic problems, the resulting hybrid scheme should be conservative in any case.
We have maintained this balance, and next, we will briefly describe our approach.
Recent approaches to solving numerical PDEs include neural network-based WENO methods that modify coefficients and smoothness indicators of established state-of-the-art numerical methods to further improve these schemes, especially near shocks. 
However, some methods achieve only first-order accuracy \cite{Stevens2020}.

In this paper, we present a new approach called "WENO-DS", a Deep learning-based extension of the family of WENO methods, and extend it to solving a general two-dimensional system of hyperbolic conservation laws
\begin{equation} \label{eq:HCL}
   U_t + F(U)_x + G(U)_y = 0.
\end{equation}
To this end, we modify the smoothness indicators of the WENO schemes using a small neural network, maintaining high accuracy in smooth regions and reducing diffusion and overshoots (oscillatory behavior) near shocks.
The resulting machine learning-enhanced WENO scheme combines accuracy and improved qualitative behavior for both smooth and discontinuous solutions.

The paper is organized as follows. 
In Section~\ref{sec:S2}, we introduce two underlying WENO schemes and explain the basic ideas, such as the smoothness indicators, on a 1D conservation law.
In Section~\ref{sec:S3}, we present our method for improving these schemes using a deep learning approach to modify the smoothness indicators accordingly.
This novel idea does not destroy the basic structure of the WENO schemes, such as the conservative property, and qualitatively improves the solution near shocks with only small additional computational costs. 
In this section, we also elaborate on implementation aspects, such as adaptive activation functions, the design of the small network, and the training procedure.
In Section~\ref{sec:S4}, we briefly describe our application example of the 
2D Euler equations of gas dynamics. 
Subsequently, in Section~\ref{sec:S5} we present in detail the numerical results with a wide range of test configurations.
Finally, in Section~\ref{sec:S6} we conclude our work and give a brief overview of future research directions.

\section{The WENO scheme} \label{sec:S2}
We first introduce the standard fifth-order WENO scheme for solving one-dimensional hyperbolic conservation laws 
\begin{equation} \label{eq:HCL_1D}
   u_t + f(u)_x = 0,
\end{equation}
as developed by Jiang and Shu \cite{jiang1996, shu1998}.
For this purpose, we consider the uniform grid defined by the points
$x_i = x_0+i\Delta x$ with cell boundaries $x_{i+\frac{1}{2}} = x_i+\frac{\Delta x}{2}$, $i = 0,\ldots,I$. 
The semi-discrete formulation of \eqref{eq:HCL_1D} can be written as
\begin{equation} \label{eq:hypsemi}
    \dt{u_i(t)} = -\frac{1}{\Delta x}\bigl(\hat{f}_{i+\frac{1}{2}}-\hat{f}_{i-\frac{1}{2}}\bigr),
\end{equation}
where $u_i(t)$ approximates $u(x_i,t)$ pointwise and $\hat{f}$ is a numerical approximation of the flux function $f$, i.e.\ $\hat{f}_{i+\frac{1}{2}}$ and $\hat{f}_{i-\frac{1}{2}}$ are numerical flux approximations at the cell boundaries $x_{i+\frac{1}{2}}$ and $x_{i-\frac{1}{2}}$, respectively.
The numerical flux $\hat{f}_{i+\frac{1}{2}}$ is chosen such that for all sufficiently smooth $u$
\begin{equation}
    \frac{1}{\Delta x} \Bigl( \hat{f}_{i+\frac{1}{2}} - \hat{f}_{i-\frac{1}{2}} \Bigr) 
    = \bigl(f(u)\bigr)_{x}\big\vert_{x=x_i} + O(\Delta x^5),
\end{equation}
with fifth-order of accuracy.
Defining a function $h$ implicitly by
\begin{equation}\label{eq:int}
    f\bigl(u(x)\bigr) = \frac{1}{\Delta x} \int_{x-\frac{\Delta x}{2}}^{x+\frac{\Delta x}{2}} h(\xi)\,d\xi,
\end{equation}
we obtain 
\begin{equation} \label{eq:approx}
    f'\bigl(u(x_i)\bigr) 
    = \frac{1}{\Delta x}\bigl(h_{i+\frac{1}{2}} - h_{i-\frac{1}{2}}\bigr),
    \qquad h_{i\pm\frac{1}{2}} = h(x_{i\pm\frac{1}{2}}),
\end{equation}
where $h_{i\pm\frac{1}{2}}$ approximates the numerical flux $\hat{f}_{\pm\frac{1}{2}}$ 
with the fifth-order of accuracy in the sense that
\begin{equation} 
    \hat{f}_{i\pm\frac{1}{2}} = h_{i\pm\frac{1}{2}} + O(\Delta x^5).
\end{equation}
This procedure results in a \textit{conservative} numerical scheme.

To ensure numerical stability, the \textit{flux splitting method} is applied. 
We therefore write the flux in the form
\begin{equation} \label{eq:fluxsplit}
    f(u) = f^+(u) + f^-(u),\quad
    \text{where}\quad \dd{f^+(u)}{u}\ge0\quad
    \text{and}\quad \dd{f^-(u)}{u}\le0.
\end{equation}
The numerical flux $\hat{f}_{i\pm\frac{1}{2}}$ is then given by 
$\hat{f}_{i\pm\frac{1}{2}} = \hat{f}_{i\pm\frac{1}{2}}^+ + \hat{f}_{i\pm\frac{1}{2}}^-$ 
and we get the final approximation
\begin{equation} \label{eq:hypsemi_splitting}
    \dt{u_i}= -\frac{1}{\Delta x}\biggl[\Bigl(\hat{f}_{i+\frac{1}{2}}^+-\hat{f}_{i-\frac{1}{2}}^+\Bigr)+\Bigl(\hat{f}_{i+\frac{1}{2}}^--\hat{f}_{i-\frac{1}{2}}^-\Bigr)\biggr].
\end{equation}
\begin{remark}
In our implementation, we use the Lax-Friedrichs flux splitting
\begin{equation} \label{eq:LF_flux_splitting}
    f^\pm(u) = \frac{1}{2}\bigl(f(u)\pm\alpha u\bigr),
\end{equation}
with $\alpha = \max\limits_u|f'(u)|$.
\end{remark}


\subsection{The fifth order WENO scheme} \label{sec:S2.1}

First, we consider the construction of  $\hat{f}^+_{i+\frac{1}{2}}$ and drop the superscript $^+$ for simplicity. 
For this approximation a 5-point stencil
\begin{equation} \label{eq:stencil_plus}
    S(i)=\{x_{i-2},\dots,x_{i+2} \}
\end{equation}
is used. 
The main idea of the fifth-order WENO scheme is to divide this stencil \eqref{eq:stencil_plus} into three candidate substencils, which are given by
\begin{equation} \label{eq:substencils}
    S^m(i)=\{x_{i+m-2}, x_{i+m-1}, x_{i+m} \}, \quad  m=0,1,2.
\end{equation}
The numerical fluxes $\hat{f}^m(x_{i+\frac{1}{2}}) = \hat{f}^m_{i+\frac{1}{2}} = h_{i+\frac{1}{2}}+ O(\Delta x^3)$ are then calculated for each of the small substencils \eqref{eq:substencils}.
Let $\hat{f}^m(x)$ be the polynomial approximation of $h(x)$ on each of the substencils \eqref{eq:substencils}. 
By evaluation of these polynomials at $x = x_{i+\frac{1}{2}}$ the following explicit formulas can be obtained \cite{shu1998}
\begin{equation} \label{eq:sub_fluxes} 
\begin{split}
   \hat{f}^0_{i + \frac{1}{2}} &= \frac{2f(u_{i-2})-7f(u_{i-1})+11f(u_{i})}{6},  \\
   \hat{f}^1_{i + \frac{1}{2}} &= \frac{-f(u_{i-1})+5f(u_{i})+2f(u_{i+1})}{6}, \\
   \hat{f}^2_{i + \frac{1}{2}} &= \frac{2f(u_{i})+5f(u_{i+1})-f(u_{i+2})}{6}, 
\end{split}
\end{equation} 
where the value of a function $f$ at $u(x_i)$ is indicated by $f(u_i)=f(u(x_i))$.
Then, we obtain a final approximation on a big stencil \eqref{eq:stencil_plus} as a linear combination of the fluxes \eqref{eq:sub_fluxes}
\begin{equation}
     \hat{f}_{i + \frac{1}{2}} = \sum_{m=0}^2 d_m  \hat{f}^m_{i+\frac{1}{2}},
\end{equation}
where the coefficients $d_m$ are the linear weights, 
which would form the upstream fifth-order central scheme for the 5-point stencil and their values are
\begin{equation} 
    d_0=\frac{1}{10}, \quad d_1=\frac{6}{10}, \quad d_2=\frac{3}{10}.
\end{equation}
As described in \cite{jiang1996, shu1998}, the linear weights can be replaced by \textit{nonlinear weights} $\omega_m^{JS}$, $m=0,1,2$, such that
\begin{equation} \label{eq:omega_flux}
    \hat{f}_{i + \frac{1}{2}} = \sum_{m=0}^2 \omega_m^{JS}\hat{f}^m_{i + \frac{1}{2}},
\end{equation}
with
\begin{equation} \label{eq:omegas}
   \omega_m^{JS} = \frac{\alpha_m^{JS}}{\sum_{i=0}^2 \alpha_i^{JS}}, 
   \quad \text{ where } \quad  
   \alpha_m^{JS} = \frac{d_m}{ (\epsilon + \beta_m)^2 }.
\end{equation}
The parameter $\beta_m$ is crucial for deciding which substencils to include in the final flux approximation. 
It is referred to as \textit{smoothness indicator} and its main role is to reduce or remove the contribution of the substencil $S^m$, which contains the discontinuity. In this case, the corresponding nonlinear weight $\omega_m^{JS}$ becomes smaller. For smooth parts of the solution, the indicators are designed to come closer to zero, so that the nonlinear weights $\omega_m^{JS}$ come closer to the ideal weights $d_m$. We will further analyze the smoothness indicators in the next section.
The parameter $\epsilon$ is used to prevent the denominator from becoming zero. In all our experiments, we set the value of $\epsilon$ to $10^{-6}$.

\subsection{Smoothness indicators}  \label{S:2.2}
In \cite{jiang1996}, the smoothness indicators have been developed as: 
\begin{equation} \label{eq:betas}
   \beta_m = \sum_{q=1}^2 \Delta x^{2q-1} \int_{x_{i-\frac{1}{2}}}^{x_{i+\frac{1}{2}}} 
      \Bigl(\ddq{q}{\hat{f}^m(x)}{x}\Bigr)^2\,dx,
\end{equation} 
with $\hat{f}^m(x)$ being the polynomial approximation in each of three substencils. 
Their explicit form corresponding to the flux approximation $\hat{f}_{i+\frac{1}{2}}$ can be obtained as 
\begin{equation}  \label{eq:betas_explicit}
\begin{split}
   \beta_0 &= \frac{13}{12} \bigl(f(u_{i-2})-2f(u_{i-1})+f(u_{i})\bigr)^2 +\frac{1}{4}\bigl(f(u_{i-2})-4f(u_{i-1})+3f(u_{i})\bigr)^2 ,  \\
   \beta_1 &= \frac{13}{12} \bigl(f(u_{i-1})-2f(u_{i})+f(u_{i+1})\bigr)^2 +\frac{1}{4}\bigl(-f(u_{i-1})+f(u_{i+1})\bigr)^2, \\
    \beta_2 &= \frac{13}{12} \bigl(f(u_{i})-2f(u_{i+1})+f(u_{i+2})\bigr)^2 +\frac{1}{4}\bigl(3f(u_{i})-4f(u_{i+1})+f(u_{i+2})\bigr)^2.
\end{split}
\end{equation} 

\begin{remark}
As mentioned before, we only considered the construction of the numerical flux $\hat{f}^+_{i+\frac{1}{2}}$. 
For the numerical approximation of the flux $\hat{f}^+_{i-\frac{1}{2}}$ we can use formulas \eqref{eq:sub_fluxes}--\eqref{eq:omegas} 
and \eqref{eq:betas_explicit} and shift each index by $-1$.
\end{remark}

The negative part of the flux splitting can be obtained using symmetry (see, e.g., \cite{wang2007}), and we briefly summarize the formulas for $\hat{f}_{i+\frac{1}{2}}^-$ and omit the superscript $^-$:
\begin{equation}  
\begin{split}
   \hat{f}^0_{i + \frac{1}{2}} &= \frac{11f(u_{i+1})-7f(u_{i+2})+2f(u_{i+3})}{6},  \\
   \hat{f}^1_{i + \frac{1}{2}} &=  \frac{2f(u_{i})+5f(u_{i+1})-f(u_{i+2})}{6}, \\
   \hat{f}^2_{i + \frac{1}{2}} &=  \frac{-f(u_{i-1})+5f(u_{i})+2f(u_{i+1})}{6}, 
\end{split}
\end{equation} 
where the weights $\omega_m^{JS}$ are computed as in \eqref{eq:omegas} using the smoothness indicators given by
\begin{equation}\label{eq:betas_explicit_minus}
\begin{split}
   \beta_0 &= \frac{13}{12} \big(f(u_{i+1})-2f(u_{i+2})+f(u_{i+3})\big)^2 +\frac{1}{4}\big(3f(u_{i+1})-4f(u_{i+2})+f(u_{i+3})\big)^2 ,  \\
   \beta_1 &= \frac{13}{12} \big(f(u_{i})-2f(u_{i+1})+f(u_{i+2})\big)^2 +\frac{1}{4}\big(f(u_{i})-f(u_{i+2})\big)^2, \\
   \beta_2 &= \frac{13}{12} \big(f(u_{i-1})-2f(u_{i})+f(u_{i+1})\big)^2 +\frac{1}{4}\big(f(u_{i-1})-4f(u_{i})+3f(u_{i+1})\big)^2.  
\end{split}
\end{equation} 
In the next section, where the deep learning algorithm will be introduced, this will help to understand how the improved smoothness indicators will be constructed.

\subsection{The WENO-Z scheme} \label{sec:S2.3}
Borges et al.\ \cite{Borges2008} pointed out that the classical WENO-JS scheme described in previous sections looses the fifth-order 
accuracy at the critical points where $f'(u)=0$, and proposed new nonlinear weights defined by
\begin{equation} \label{eq:WENOZ}
   \omega_m^Z =  \frac{\alpha^Z_m}{\sum\limits_{i=0}^2\alpha^Z_i}, \quad \text{ where } \quad \alpha^Z_m
    = d_m \biggl[ 1+ \Bigl(\frac{\tau_5}{\beta_m + \epsilon} \Bigr)^2 \biggr]
\end{equation}
and
\begin{equation} \label{eq:tau}
    \tau_5 = |\beta_0 - \beta_2|
\end{equation}
is a new global smoothness indicator. 

\section{Deep smoothness WENO scheme} \label{sec:S3}
In \cite{kossaczka2021, kossaczka2022, kossaczka2022deep} the new WENO-DS scheme based on the improvement of the smoothness indicators was developed.
The smoothness indicators $\beta_m$, $m=0,1,2$, are multiplied by the perturbations $\delta_m$, 
which are the outputs of the respective neural network algorithm. 
The new smoothness indicators are denoted by $\beta_m^{DS}$:
\begin{equation} \label{eq:betas_DS}
    \beta_m^{DS} = \beta_m (\delta_m + C), \qquad m=0,1,2,
\end{equation}
where $C$ is a constant that ensures the consistency and accuracy of the new method. 
In all our experiments we set $C=0.1$.
For more details and corresponding theoretical proofs of accuracy and consistency, we refer to \cite{kossaczka2021, kossaczka2022}.

Note that the formulation of the new smoothness indicators $\beta_m^{DS}$ as a multiplication of the original ones with the perturbations $\delta_m$ is very favorable. 
In a case where the original smoothness indicator converges to zero, the improved $\beta_m^{DS}$ behaves in the same way. 
On the other hand, if a subset $S^m$ contains a discontinuity, the perturbation $\delta_m$ can improve the original smoothness indicator so that the final scheme exhibits better numerical approximations.
Moreover, the theoretical convergence properties are not lost, see \cite{kossaczka2021, kossaczka2022}.

In \cite{kossaczka2021} the algorithm was successfully applied to one-dimensional benchmark examples such as the Burgers' equation, 
the Buckley-Leverett equation, the one-dimensional Euler system, and the two-dimensional Burgers' equation. 
In \cite{kossaczka2022}, the algorithm was extended to nonlinear degenerate parabolic equations and further applied to computational finance problems in \cite{kossaczka2022deep}. 
The theoretical order of convergence was demonstrated on the smooth solutions and the large numerical improvements were obtained when comparing the WENO-DS method with the original WENO methods.

\subsection{Preservation of a conservative property for WENO-DS scheme} \label{sec:S3.1}
However, the multipliers introduced for the smoothness indicators in \cite{kossaczka2021} were cell-based (not interface-based). 
This means that although the high numerical accuracy was theoretically demonstrated and numerically confirmed, the guarantee of the conservative property was lost. 
As stated in \cite{kossaczka2022}, the conservative property can be easily recovered by defining the multipliers such that
\begin{equation} \label{eq:betas_DS_cons}
\begin{aligned} 
     \beta_{m,i+\frac{1}{2}}^{DS} &= \beta_{m,i+\frac{1}{2}} (\delta_{m,i+\frac{1}{2}} + C), \\
          \beta_{m,i-\frac{1}{2}}^{DS} &= \beta_{m,i-\frac{1}{2}} (\delta_{m,i-\frac{1}{2}} + C), 
\end{aligned}
\end{equation}
with
\begin{equation} \label{eq:deltas}
    \delta_{0,i+\frac{3}{2}} =  \delta_{1,i+\frac{1}{2}} =  \delta_{2,i-\frac{1}{2}}, \quad i=0, \ldots, N.
\end{equation}
This makes the multipliers depend on the location of the substencils corresponding to $\beta_{m,i+\frac{1}{2}}$ and $\beta_{m,i-\frac{1}{2}}$. 
This ensures that the values $\hat{f}^\pm_{i-\frac{1}{2}}$ can be obtained from the values $\hat{f}^\pm_{i+\frac{1}{2}}$ by simple index shifting and that the conservative property is preserved.

\subsection{Structure of neural network} \label{sec:S3.2}
To ensure the consistency of a numerical method, the Convolutional Neural Network (CNN) is used. 
This is crucial to ensure the spatial invariance of the resulting numerical method.  
This means that the multipliers $\delta_m$ are independent of their position in the spatial grid and only depend on the solution itself.

\begin{figure}[ht!] 
\centering\includegraphics[width=1\linewidth]{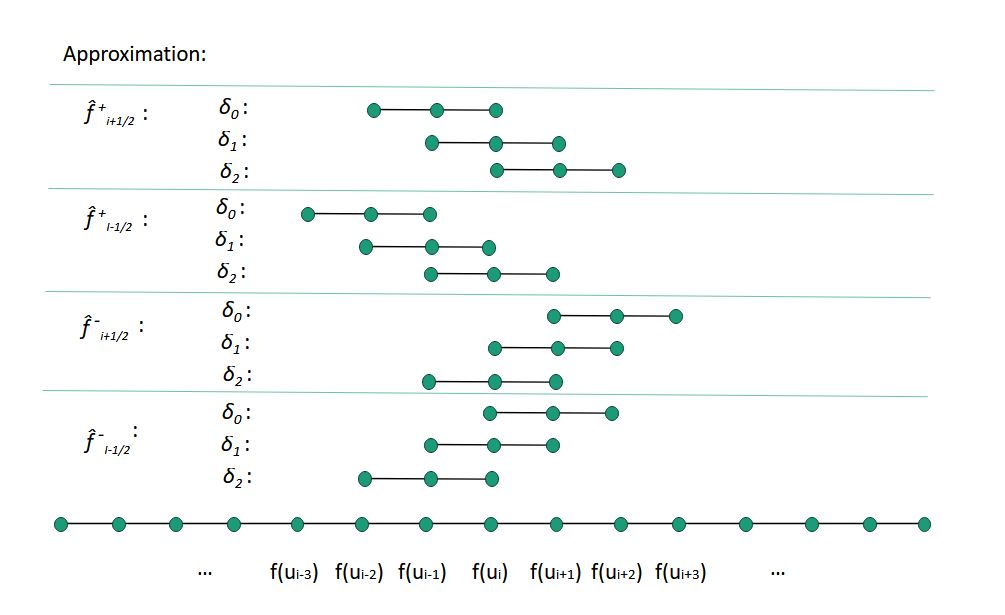}
\caption{The substencils used for computation of multipliers $\delta_m$, $m=0,1,2$ corresponding to the flux approximations $\hat{f}^\pm_{i\pm\frac{1}{2}}$, assuming that for the receptive field of the CNN holds $2k+1=3$.}
\label{fig:Deltas_stencils} 
\end{figure}

Let us formulate the CNN as a function $H(\cdot)\in\mathbb{R}^{2k+1} \to \mathbb{R}$, where $2k+1$ is the size of the receptive field of the CNN:
\begin{equation} \label{eq:CNN_func}
H\bigl(\bar{f}(\bar{u}_{i})\bigr) = {\rm CNN} \bigl(\bar{f}(\bar{u}_{i})\bigr).
\end{equation}
As an input we define a vector
\begin{equation} \label{eq:cnn_input}
\begin{split}
\bar{f}(\bar{u}_i) &= \bigl(f(u(x_{i-k})), f(u(x_{i-k+1})), \ldots, f(u(x_{i+k}))\bigr), \\
\bar{u}_i &= \bar{u}(\bar{x}_i)=\bigl(u(x_{i-k}), u(x_{i-k+1}), \ldots, u(x_{i+k})\bigr).
\end{split}
\end{equation}
The Figure~\ref{fig:Deltas_stencils} shows the values from which the multipliers $\delta_m$, $m=0,1,2$ are constructed, assuming $2k+1=3$ for the receptive field. 
In this case, the values used to compute the original smoothness indicators are also used to compute the multipliers $\delta_m$, $m=0,1,2$, (see equations \eqref{eq:betas_explicit} and \eqref{eq:betas_explicit_minus}).
If we enlarge the receptive field of the CNN, we also enlarge the stencil for computing the multipliers $\delta_m$, $m=0,1,2$.
In this way, the smoothness indicators are basically computed from a wider stencil, which can lead to better numerical approximations.
In this case, we just need to add more bounds before feeding the values \eqref{eq:cnn_input} to the CNN.
 
As we are improving the existing numerical scheme and adding a neural network part to it, 
it is important that the new numerical scheme remains computationally efficient. 
The neural network part added to the numerical scheme could be computationally expensive. 
However, we propose to use only a small CNN, which would not have such high computational costs. 
The detailed structure of the CNN can be found in Section~\ref{sec:S4.1}.

It was pointed out in \cite{kossaczka2021} that better numerical results were obtained using two different neural networks for the positive and negative part of a flux.
We experimentally found that we can avoid using more neural networks and use only one CNN.
On the other hand, we can achieve better results by using a superior training procedure and adaptive activation functions.
More details will be discussed in the next subsections.

For convergence and consistency of the numerical scheme, all hidden layers of the CNN must be differentiable functions, and the activation function in the last CNN layer must be bounded from below \cite{kossaczka2022}.
Experimentally, we found that the use of a \textit{softplus} activation function in the last CNN layer is more effective and gives better numerical results compared to e.g.\ \textit{sigmoid} as used in \cite{kossaczka2021}.

\subsubsection{Adaptive activation functions}
We can make the training more effective and get better numerical results by using \textit{adaptive activation functions}. \cite{jagtap2020adaptive,jagtap2020locally,jagtap2022deepR}. 
The activation function is one of the most important hyperparameters in neural network architectures. 
The purpose of this hyperparameter is to introduce nonlinearity into the prediction. 
There are many activation functions proposed in the literature; see the comprehensive survey \cite{jagtap2023important} for more details. 
However, there is no basic rule for the choice of the activation function. 
This is the motivation to use an adaptive activation function that can adapt to the problem at hand. 
In this work, we used global adaptive activation functions \cite{jagtap2020adaptive}, where the additional slope parameter is introduced in the activation function as follows.

For the ELU activation function, we train the additional parameter $\alpha$:
\begin{equation}
   \text{ELU} =  \begin{cases}
            x, \qquad &\text{if} \quad x>0, \\
            \alpha (\exp(x)-1) \qquad &\text{if} \quad x \leq 0
   \end{cases}
\end{equation}
and we denote the adaptive ELU as \textit{aELU}.
For the softplus activation function, we train the additional parameter $\beta$:
\begin{equation}
   \text{Softplus}(x) = \frac{1}{\beta} \log(1 + \exp(\beta x)) 
\end{equation}
and we denote the adaptive softplus as \textit{aSoftplus}.


\subsection{Training procedure} \label{sec:S3.3}

In this section, we describe how the training procedure for WENO-DS is carried out.
We have experimented with different training procedures and have found experimentally that following the training procedure described in \cite{kossaczka2022} gives the best numerical results.
First, we have to create the data set. 
For this purpose we compute the reference solutions using the WENO-Z method on a fine grid of $I\times J = 400\times400$ space points up to the given final time $T$, where $t_n$ represents the time points, $n=0,\ldots,N$. 
More details on the construction of the reference solutions are given in Sections~\ref{sec:S5.1}, \ref{sec:S5.2}, \ref{sec:S5.3}.  

During training, we compute the numerical solutions on a grid of $I\times J = 100\times100$ space points.
At the beginning of a training we randomly select a problem from a data set and perform a single time step to get to the time $t_{n+1}$, using CNN to predict the multipliers $\delta_m$. 
However, by performing a single time step on a coarse grid, we do not match the time step size of the fine precomputed solutions, 
as the adaptive time step size is used. 
So we simply take the closest reference solution from the data set, use it as an initial condition, and do another small time step to get a reference solution in time $t_{n+1}$. 
Then we compute the loss and its gradient with respect to the weights of the CNN. 

We then decide whether to proceed to the next time step of a current problem or to choose another problem from our dataset and run a time step of that problem. 
The probability of choosing the new problem has to be determined at the beginning of the training session and we set it to $\varphi = 0.5$ in our experiments.
We set the maximum number of opened problems to $150$. 
We remember all opened problems, and if no new problem is opened (with probability $1-\varphi$), 
or if the maximum number of opened problems is reached, we execute the next time step of a problem uniformly chosen from the set of already opened problems. 
 Keeping the solution from the previous time step as initial data, we repeat the same procedure until we reach the maximum number of training steps.
This training procedure gives us a great opportunity to mix the solutions with different initial data and in different time points, which makes the training more effective.

\subsubsection{Optimizer and the optimal learning rate}
To train the network, we used a gradient-based optimizer, namely a variant of stochastic gradient descent, the Adam optimizer \cite{kingma2014adam}.

The learning rate is another important hyperparameter to choose. 
A larger learning rate may miss the local minima, and a smaller learning rate may require a large number of iterations to reach convergence. 
Therefore, it is important to find a near-optimal learning rate.
In this work, the learning rate is $0.001$ to update the weights of the CNN. 
This near-optimal learning rate was found through experiments.

\subsubsection{Loss function}
In this work, the loss function consists of the data mismatch term between the solution predicted by the networks and the reference solution.
For the loss function, we use the mean square error loss as follows:
\begin{equation} \label{eq:loss_L2}
   LOSS_{\rm MSE}(u) = \frac{1}{I} \sum_{i=0}^I (u_i - u_i^{\rm ref})^2,
\end{equation}
where $u_i$ is a numerical approximation of $u(x_i)$ and $u_i^{\rm ref}$ is the corresponding reference solution. 
The $L_2$ norm-based loss function has the advantage of stronger gradients with respect to $u_i$, resulting in faster training.
However, in our examples, we use the $L_1$ norm as the main error measure, 
which is more typical for measuring errors for hyperbolic conservation laws. 
Thus, for validation during training, we use the metrics
\begin{equation} \label{eq:loss_L1}
   L_1(u) = \frac{1}{I} \sum_{i=0}^I |u_i - u_i^{\rm ref}|.
\end{equation}

\section{Application of our approach to the 2D Euler equations} \label{sec:S4}
We consider the two-dimensional Euler equations of gas dynamics in the form \eqref{eq:HCL} with
\begin{equation} \label{eq:Euler_HCL}
    U=\begin{pmatrix}
    \rho \\ \rho u \\\rho v \\ E
\end{pmatrix} \qquad
    F(U)=\begin{pmatrix}
    \rho u \\ \rho u^2 + p \\\rho u v \\ u (E +  p)
\end{pmatrix} \qquad
    G(U)=\begin{pmatrix}
    \rho v \\\rho u v \\ \rho v^2 + p \\ v (E +  p)
\end{pmatrix} \qquad
\end{equation}
for polytropic gas.
Here, the variable $\rho$ is the density, $u$ the $x$-velocity component, $v$ the $y$-velocity component, $E$ the total energy and $p$ the pressure.
Further, it holds
\begin{equation}
        p = (\gamma - 1)\Bigl[E-\frac{\rho}{2}(u^2+v^2)\Bigr].
\end{equation}
$\gamma$ denotes the ratio of the specific heats and we will use $\gamma \in (1.1, 1.67)$ in this paper.

We consider the spatial domain $[0,1]\times[0,1]$ and solve the Riemann problem with the following initial condition
\begin{equation}  \label{eq:Riemann_IC}
  (\rho,u, v,p) = \begin{cases} 
    (\rho_1,u_1, v_1,p_1) \quad x > 0.5 \quad  \text{and} \quad  y > 0.5,\\
     (\rho_2,u_2, v_2,p_2) \quad x < 0.5 \quad  \text{and} \quad  y > 0.5,\\
      (\rho_3,u_3, v_3,p_3) \quad x < 0.5 \quad  \text{and} \quad  y < 0.5,\\
       (\rho_4,u_4, v_4,p_4) \quad x > 0.5 \quad  \text{and} \quad  y < 0.5.\\
\end{cases} 
\end{equation}

The combination of four elementary planar waves is used to define the classification of the Riemann problem. 
A detailed study of these configurations has been done in \cite{schulz1993classification, schulz1993, chang1995, chang1999, kurganov2002, zhang1990} and there are 19 different possible configurations for polytropic gas.
These are defined by three types of elementary waves, namely a backward rarefaction wave $\overleftarrow{R}$, a backward shock wave $\overleftarrow{S}$, a forward rarefaction wave $\overrightarrow{R}$, a forward shock wave $\overrightarrow{S}$ and a contact discontinuity $J^{\pm}$, where the superscript $\pm$ refers to negative and positive contacts.

To obtain the WENO approximations in the two-dimensional example, we apply the procedure described in Section~\ref{sec:S2} using the dimension-by-dimension principle. 
Thus we obtain the flux approximations for \eqref{eq:HCL} as 
\begin{equation}
\begin{split}
    \frac{1}{\Delta x} \bigl(\hat{f}_{i+\frac{1}{2}} - \hat{f}_{i-\frac{1}{2}} \bigr) 
          &= \bigl(F(U)\bigr)_{x}\big\vert_{(x_i,y_j)} + O\bigl(\Delta x^5\bigr), \\
    \frac{1}{\Delta y} \bigl( \hat{g}_{i+\frac{1}{2}} - \hat{g}_{i-\frac{1}{2}} \bigr)
           &= \bigl(G(U)\bigr)_{y}\big\vert_{(x_i,y_j)} + O\bigl(\Delta y^5\bigr),
\end{split}
\end{equation}
with the uniform grid defined by the nodes $(x_i,y_j)$, $\Delta x = x_{i+1}-x_i$, 
$\Delta y = y_{j+1}-y_j$, $i = 0,\ldots,I$, $j = 0,\ldots,J$.

In our examples, we proceed with the implementation of the Euler system using characteristic decomposition. 
This means that we first project the solution and the flux onto the characteristic fields using left eigenvectors. 
Then we apply the Lax-Friedrichs flux splitting \eqref{eq:LF_flux_splitting} for each component of the characteristic variables. 
These values are fed into the CNN and the enhanced smoothness indicators are computed.
After obtaining the final WENO approximation, the projection back to physical space is done using the right eigenvectors, see \cite{shu1992} for more details on this procedure.

\subsection{Size of the neural network} \label{sec:S4.1}

\begin{figure}[htb] 
\centering
    \begin{subfigure}{0.95\textwidth}
            \includegraphics[width=\textwidth]{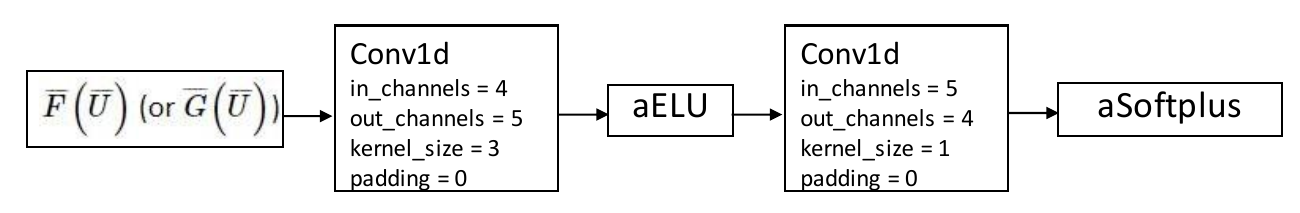}
            \caption{Two hidden layers, lower number of channels, receptive field of  size $3$.}
            \label{fig:structure_1}
    \end{subfigure}
    \begin{subfigure}{0.95\textwidth}
            \includegraphics[width=\textwidth]{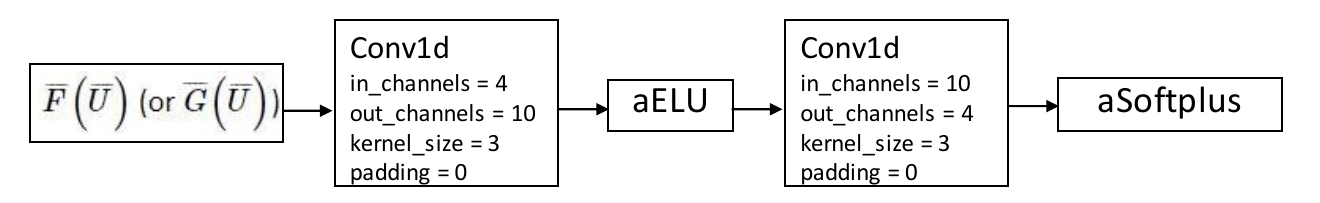}
            \caption{Two hidden layers, higher number of channels, receptive field of  size $5$.}
            \label{fig:structure_2}
    \end{subfigure}
        \begin{subfigure}{0.99\textwidth}
            \includegraphics[width=\textwidth]{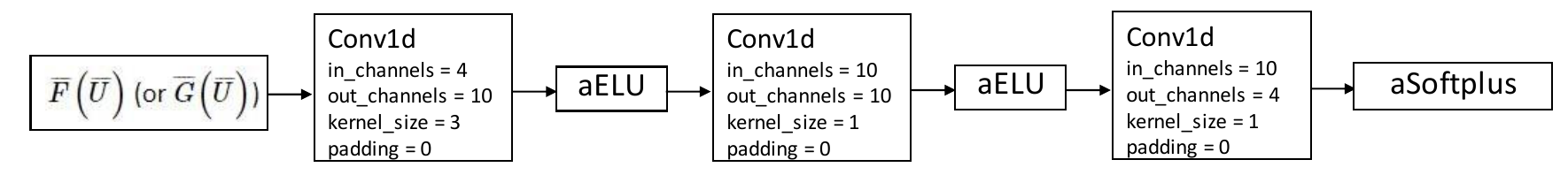}
            \caption{Three hidden layers, higher number of channels, receptive field of  size $3$.}
            \label{fig:structure_3}
    \end{subfigure}
\caption{Different structures of the convolutional neural network.}
\label{fig:structures}
\end{figure}

In our paper, we considered different structures of neural networks and carried out numerous experiments with them.
First, we used a rather simple CNN with only two layers and a receptive field of width $3$. 
The structure is shown in Figure~\ref{fig:structure_1}.  
The advantage of this is its computational efficiency. 
Second, we used a CNN with the same number of layers, but we increased the number of channels and made the receptive field wider. 
The structure is shown in Figure~\ref{fig:structure_2}. 
Finally, we used only a receptive field of width $3$, but added one more layer and used a more complex neural network, as shown in Figure~\ref{fig:structure_3}.
Each of these neural networks gave interesting results and we summarize them in Section~\ref{sec:S5}.

As can be seen, we have $4$ input channels in the first hidden layer and $4$ output channels in the last hidden layer in each CNN. 
These represent the dimension of the solution $U$ from \eqref{eq:Euler_HCL}. 
In this way, the neural network also takes in information from other variables, which can be useful for improving the numerical solution.
The input $\Bar{F}(\Bar{U})$, respectively $\Bar{G}(\Bar{U})$ represents the numerical approximation after the projection using left eigenvectors and after applying the flux splitting method.

We also have to adapt the loss function from \eqref{eq:loss_L2} and use it for training
\begin{equation} \label{eq:Euler_loss}
\begin{split}
    LOSS_{\rm MSE}(\rho, u, v, p) 
    &= LOSS_{\rm MSE}(\rho) + LOSS_{\rm MSE}(u) + LOSS_{\rm MSE}(v) + LOSS_{\rm MSE}(p)
\end{split}
\end{equation}
and for the validation during training from \eqref{eq:loss_L1} 
\begin{equation} \label{eq:Euler_validation}
      L_1(\rho, u, v, p) = L_1(\rho) + L_1(u) + L_1(v) + L_1(p).
\end{equation}
When we plot the error on validation problems, we rescale the values for each validation problem to be in the interval $[0,1]$ using the relationship
\begin{equation} \label{eq:Euler_validation_adjusted}
    L^*_1(\rho, u, v, p) = \frac{L^l_1(\rho, u, v, p)}{\max_l(L^l_1(\rho, u, v, p))}, \qquad l=0,\dots,L,
\end{equation}
where $L$ denotes the total number of training steps.

\subsection{Construction of the data set for the CNN training procedure}\label{sec:S4.2}
For each of the 19 configurations of the Riemann problem, the specific relations must be satisfied by the initial data and the symmetry properties of the solution. 
We present the formulas given in \cite{schulz1993} and create the data sets for the CNN training according to these formulas.

We define
\begin{equation}\label{eq:relations_phi_psi}
    \Phi_{lr} := \frac{2\sqrt{\gamma}}{\gamma - 1} \Big(\sqrt{\frac{p_l}{\rho_l}} - \sqrt{\frac{p_r}{\rho_r}}\Big), \quad \Psi_{lr}^2 := \frac{(p_l-p_r)(\rho_l-\rho_r)}{\rho_l \rho_r}, \quad (\Psi_{lr}>0)
\end{equation}
and
\begin{equation}\label{eq:relations_pi}
    \Pi_{lr} := \Big( \frac{p_l}{p_r} + \frac{(\gamma-1)}{(\gamma +1)}\Big) \Big/ \Big(1+\frac{(\gamma-1)}{(\gamma+1)}\frac{p_l}{p_r}\Big).
\end{equation}

In Sections \ref{sec:S5.1}, \ref{sec:S5.2}, \ref{sec:S5.3} we list the specific relations for given examples that are sufficient to uniquely define the solution. 
Following these relations, we randomly generate the initial data and construct our data sets.

\section{Numerical results} \label{sec:S5}
To demonstrate the efficiency of the proposed method, in this section, we present the numerical results obtained with the WENO-DS method after the CNN training procedure. 
Note that the CNN training procedure only needs to be performed once as \textit{offline} training for each of the examples presented in Sections~\ref{sec:S5.1}, \ref{sec:S5.2}, \ref{sec:S5.3}.
No additional training was performed for the examples in Section~\ref{sec:S5.4} as we show the results using the same trained CNN from the previous examples.
In Section~\ref{sec:S5.5} we perform two more trainings with larger CNN and illustrate the results. 
Further details can be found in the respective sections.

For the following system of ordinary differential equations (ODEs) 
\begin{equation}
        \dt{U(t)}= L(U),
\end{equation}
we use a third-order \textit{total variation diminishing} (TVD) Runge-Kutta method \cite{jiang1996} given by 
\begin{equation} \label{eq:runge_kutta}
\begin{split}
  U^{(1)} &= U^n + \Delta t \,L(U^n),  \\
  U^{(2)} &= \frac{3}{4}U^n + \frac{1}{4}U^{(1)} + \frac{1}{4}\Delta t\, L(U^{(1)}), \\
  U^{n+1} &= \frac{1}{3}U^n + \frac{2}{3}U^{(2)} + \frac{2}{3}\Delta t\,L(U^{(2)}), 
\end{split}\end{equation}
where $U^n$ is the numerical solution at the time step $n$. 

For the scheme \eqref{eq:runge_kutta} we use an adaptive step size
\begin{equation}
    \Delta t = 0.6 \min\Big(\frac{\Delta x}{a}, \frac{\Delta y}{a}\Big), 
\end{equation}
with 
\begin{equation}
    a = \max_{\substack{i=0,\ldots,I\\ j=0,\ldots,J}}(|\lambda^+_{i,j}|, |\lambda^-_{i,j}|)
    \quad \lambda^\pm = V \pm c, \quad V = \sqrt{u^2+v^2} \quad c^2= \gamma\,\frac{p}{\rho},
\end{equation}
where $u$, $v$ are the velocities and $c$ is the local speed of sound.

In the sequel we enumerate the different configurations of initial conditions according to \cite{kurganov2002}.

\subsection{Configuration 2} \label{sec:S5.1}
This is the configuration with four rarefaction waves: $\overrightarrow{R}_{21}$, $\overleftarrow{R}_{32}$, $\overleftarrow{R}_{34}$, $\overrightarrow{R}_{41}$.
The detailed analysis was done in \cite{zhang1990, schulz1993} and we have to satisfy the following relations for this case:
\begin{equation}  \label{eq:relations_R2}
\begin{split}
    &u_2-u_1 = \Phi_{21}, \quad u_4-u_3 = \Phi_{34}, \quad u_3=u_2, \quad u_4=u_1, \\
    &v_4-v_1 = \Phi_{41}, \quad v_2-v_3 = \Phi_{32}, \quad v_2=v_1, \quad v_3=v_4
\end{split}
\end{equation}
with the compatibility conditions $\Phi_{21} = -\Phi_{34}$ and $\Phi_{41}=-\Phi_{32}$.
Moreover, for a polytropic gas the equations
\begin{equation} \label{eq:relations_R2_poly}
    \rho_l/\rho_r = (p_l/p_r)^{1/\gamma} \quad \text{for} \quad (l,r) \in \{(2,1),(3,4),(3,2),(4,1)\}
\end{equation}
have to be included.
Furthermore, we have
$\rho_2=\rho_4$, $\rho_1 = \rho_3$, $p_1=p_3$, $p_2=p_4$, $u_2-u_1=v_4-v_1$ and $u_4-u_3=v_2-v_3$. 

We use for creating of the data set the values
\begin{equation} \label{eq:par_range_R2}
\begin{split}
      \rho_1 \in \mathcal{U}[0.7,2], \quad \rho_2 &\in \mathcal{U}[0.5,\rho_1], \quad p_1 \in \mathcal{U}[0.2,1.5], \\
  \quad u_1 \in \mathcal{U}[-1,1], \quad &v_1 = u_1, \quad \gamma \in (1.1, 1.67)
\end{split}
\end{equation}
and for the other values we use the relations \eqref{eq:relations_R2}, \eqref{eq:relations_R2_poly} with \eqref{eq:relations_phi_psi}. 
We also compute the reference solutions using the WENO-Z method on a grid $I\times J = 400\times400$ space points up to the final time $T \in \mathcal{U}[0.1,0.2]$ and create the data set consisting of $50$ reference solutions.

\begin{figure}[htb] 
\centering
    \includegraphics[width=0.7\textwidth]{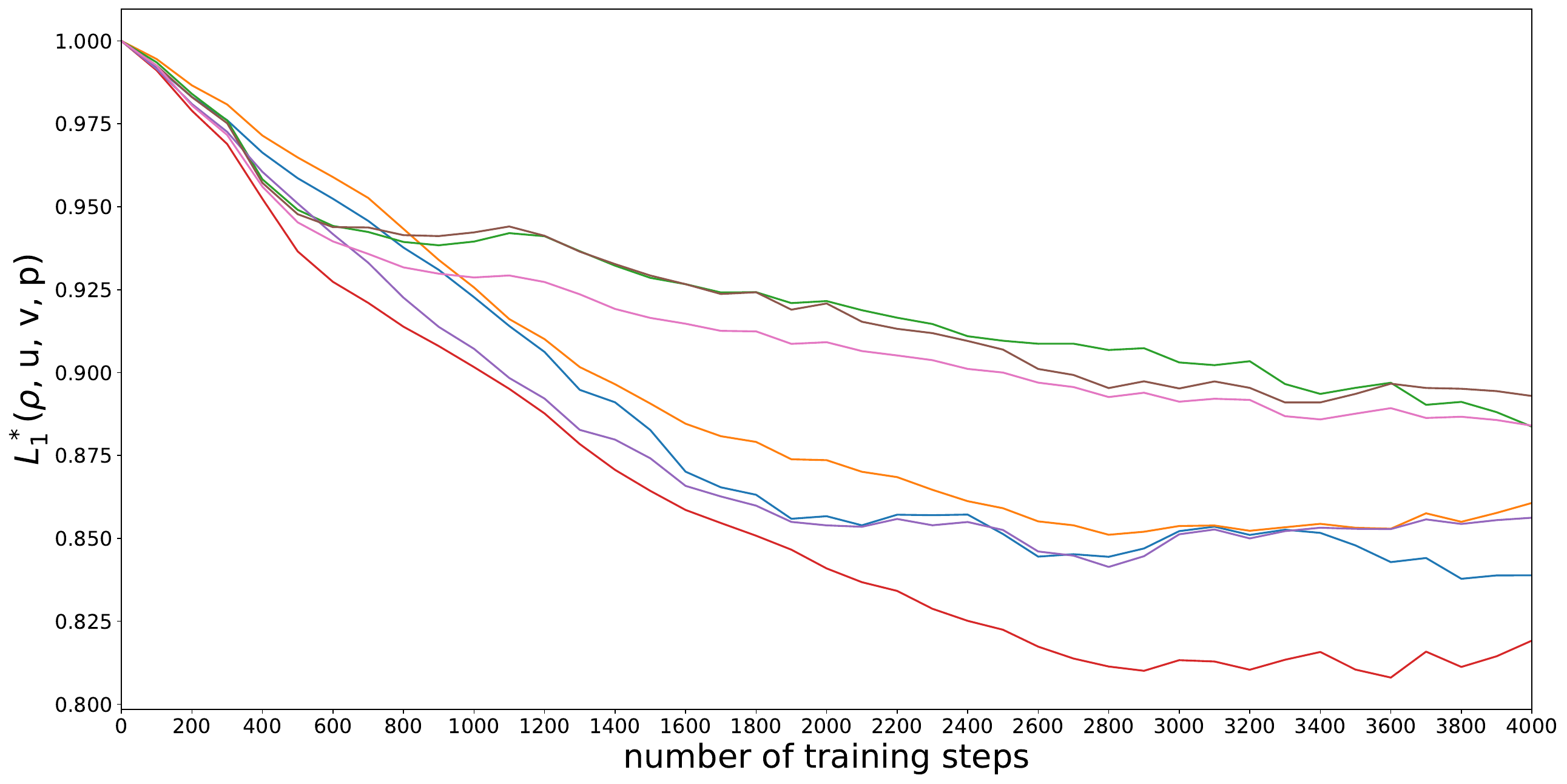}
     \caption{The values \eqref{eq:Euler_validation_adjusted} for different validation problems evaluated each $100$ training steps.}
\label{fig:validation_Riemann_2}
\end{figure}

For training, we use the training procedure described in Section~\ref{sec:S3.3}. 
First, we use the simplest neural network structure shown in Figure~\ref{fig:structure_1} and perform the training for the total number of $4000$ training steps. 
We plot the evolution of the $L^*_1$ error \eqref{eq:Euler_validation_adjusted} for the validation problems in Figure~\ref{fig:validation_Riemann_2}. 
Note that these problems were not included in the training data, and the initial conditions of these problems were generated analogously to the construction of the training data set. 
For these problems, we measured the error every $100$ training steps and at a randomly chosen final time $T$. 
We select the final model based on the evolution of the error of the validation set.
We see that the error decreases up to a certain point for all problems and then starts to increase for some problems. 
Longer training would lead to overfitting of the training data.
Finally, we choose the final model from the $2800$ training step and present the results using this model.  

\begin{figure}[h]
\begin{subfigure}{.32\textwidth}
\centering
\includegraphics[width=\linewidth]{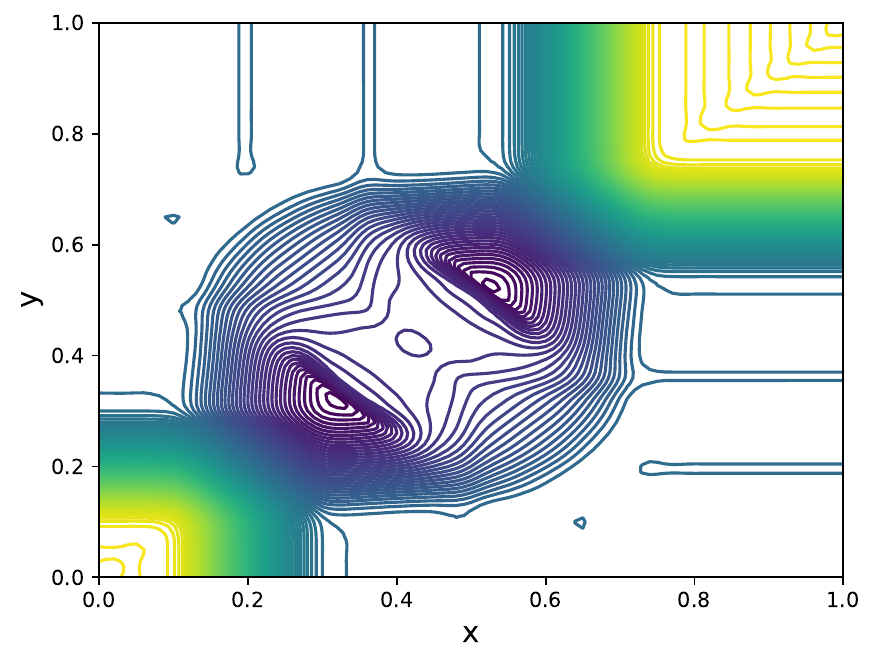}
\caption{WENO-DS}\label{fig:R2_WENODS}
\end{subfigure}
\begin{subfigure}{.32\textwidth}
\centering
\includegraphics[width=\linewidth]{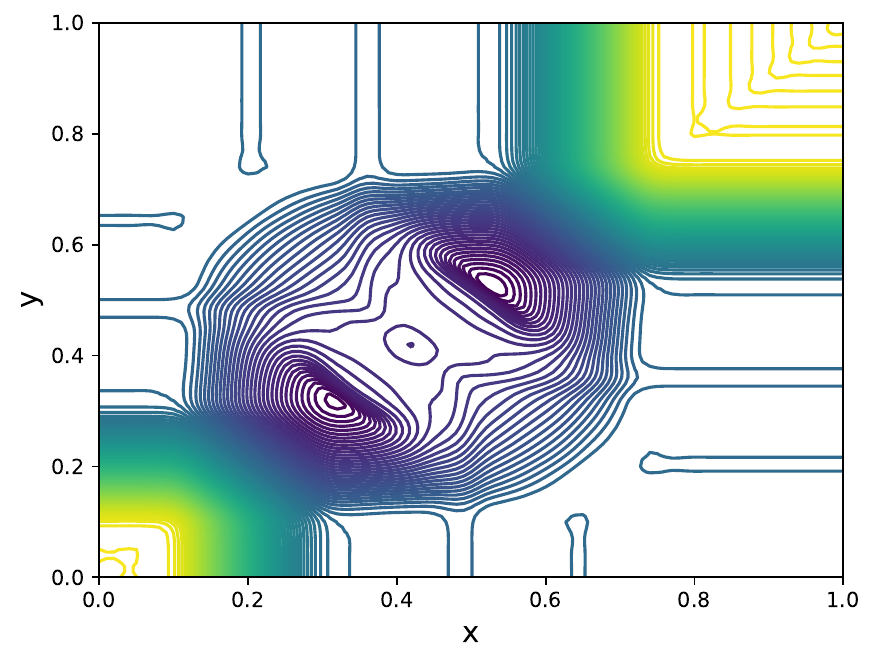} 
\caption{WENO-Z}\label{fig:R2_WENOZ}
\end{subfigure}
\begin{subfigure}{.32\textwidth}
\centering
\includegraphics[width=\linewidth]{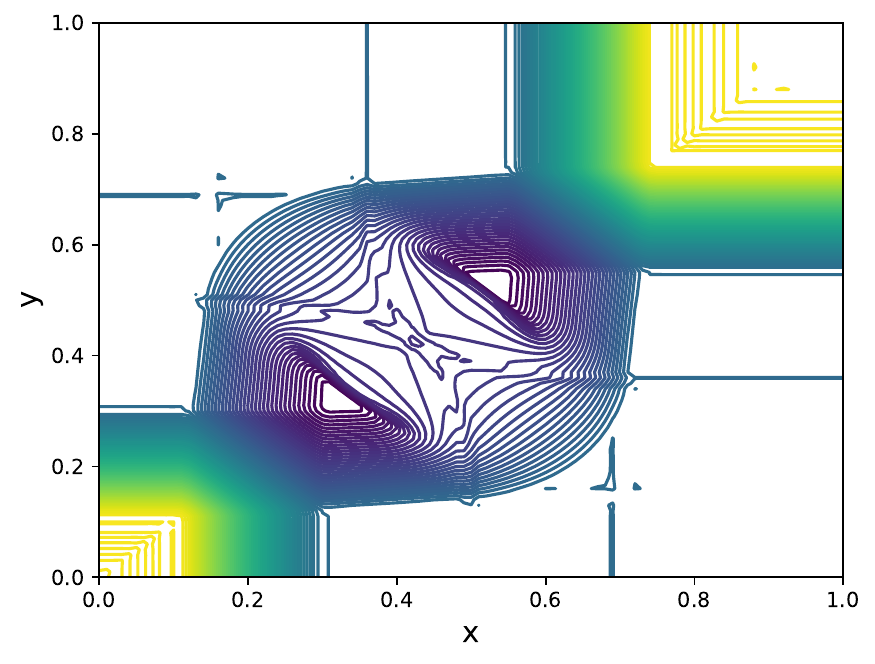} 
\caption{reference solution}\label{fig:R2_reference}
\end{subfigure}
\caption{Density contour plot for the solution of the Riemann problem with the initial condition \eqref{eq:Riemann_IC_2}, $I \times J = 100 \times 100$, $T=0.2$.} \label{fig:Riemann_2}
\end{figure}

\begin{figure}[h!]
\centering
\begin{subfigure}{.4\textwidth}
\centering
\includegraphics[width=\linewidth]{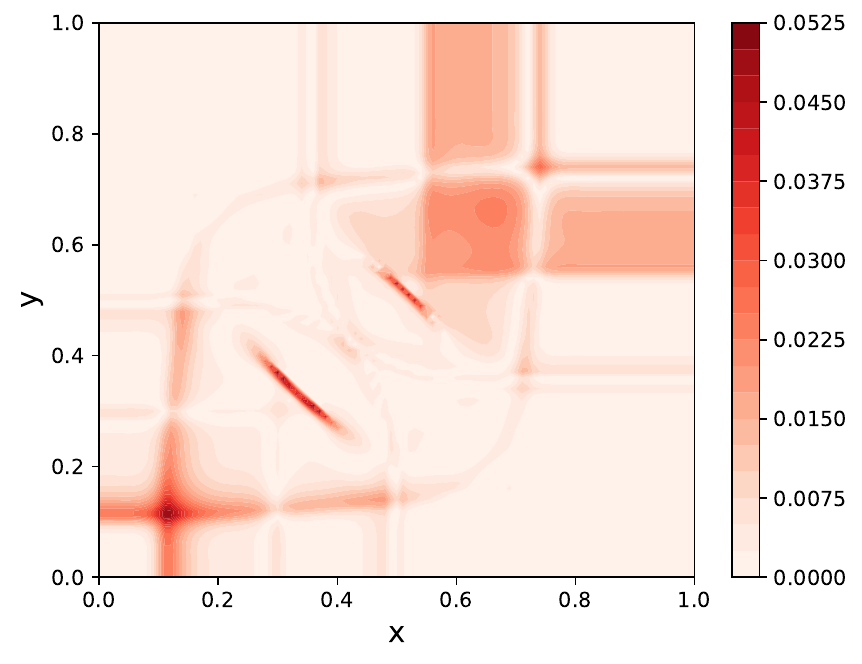}
\caption{WENO-DS}\label{fig:R2_WENODS_error}
\end{subfigure}
\begin{subfigure}{.4\textwidth}
\centering
\includegraphics[width=\linewidth]{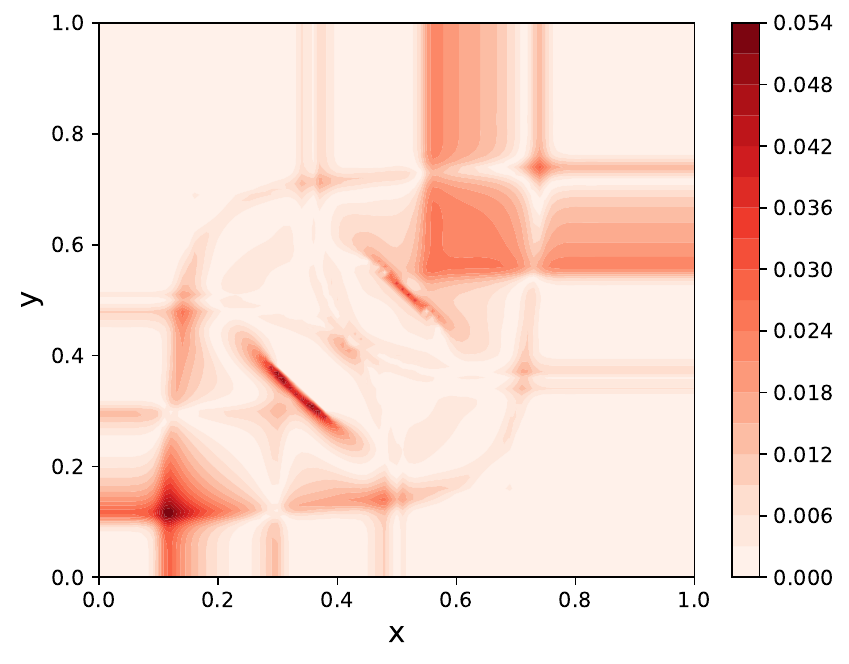} 
\caption{WENO-Z}\label{fig:R2_WENOZ_error}
\end{subfigure}
\caption{Absolute pointwise errors for the density solution of the Riemann problem with the initial condition \eqref{eq:Riemann_IC_2}, $I \times J = 100 \times 100$, $T=0.2$.} \label{fig:Riemann_2_error}
\end{figure}

As a test problem we use the problem from \cite{kurganov2002} with $\gamma = 1.4$, $T=0.2$ and the initial condition
\begin{equation}  \label{eq:Riemann_IC_2}
  (\rho,u, v,p) = \begin{cases} 
    (1, 0, 0, 1) \quad &x > 0.5 \quad  \text{and} \quad  y > 0.5,\\
     (0.5197,-0.7259, 0, 0.4) \quad &x < 0.5 \quad  \text{and} \quad  y > 0.5,\\
      (1, -0.7259, -0.7259, 1) \quad &x < 0.5 \quad  \text{and} \quad  y < 0.5,\\
       (0.5197, 0, -0.7259, 0.4) \quad &x > 0.5 \quad  \text{and} \quad  y < 0.5.\\
\end{cases} 
\end{equation}
The results are shown in Table~\ref{tab:Riemann_2}. 
As can be seen, we achieve a significant error improvement for all four variables and for different discretizations. 
It should be noted that we trained only with the discretization $100 \times 100$ space points and did not retrain the neural network for different discretizations. 
We refer to the error of the WENO-Z method divided by the error of WENO-DS (rounded to 2 decimal points) as the 'ratio'. 
The density contour plots are shown in Figure~\ref{fig:Riemann_2} and the absolute pointwise errors for the density are shown in Figure~\ref{fig:Riemann_2_error}.

\begin{table}[h] 
 \centering
    \scalebox{0.7}{ \begin{tabular}{|c|c|c|c|c|c|c|c|c|c|}
    \hline
    \multicolumn{1}{|c|}{$I \times J$}&\multicolumn{3}{|c|}{\ $50 \times 50$}&\multicolumn{3}{|c|}{\ $100 \times 100$} &\multicolumn{3}{|c|}{\ $200 \times 200$}\\
    \hline
        & WENO-Z &  WENO-DS & ratio & WENO-Z &  WENO-DS & ratio& WENO-Z &  WENO-DS & ratio\\
    \hline
    \hline
 $\rho$ &  0.012488 &  0.010722 &  1.16   &  0.005465 &  0.004648 &  1.18 &  0.001862 &  0.001547 &  1.20  \\ \hline
$u$ &  0.014363 &  0.011986 &  1.20  &  0.006153 &  0.005066 &  1.21 &  0.002053 &  0.001627 &  1.26   \\ \hline
$v$ &  0.014363 &  0.011986 &  1.20   &  0.006153 &  0.005066 &  1.21  &  0.002053 &  0.001627 &  1.26  \\ \hline
 $p$ &  0.013113 &  0.011510 &  1.14  &  0.005619 &  0.004899 &  1.15 &  0.001879 &  0.001587 &  1.18   \\ \hline
\end{tabular} }
    \caption{Comparison of $L_1$ error of WENO-Z and WENO-DS methods for the solution of the Euler system with the initial condition \eqref{eq:Riemann_IC_2} for different spatial discretizations, $T=0.2$.}
\label{tab:Riemann_2}
\end{table}

\begin{figure}[htb] 
\centering
    \includegraphics[width=0.5\textwidth]{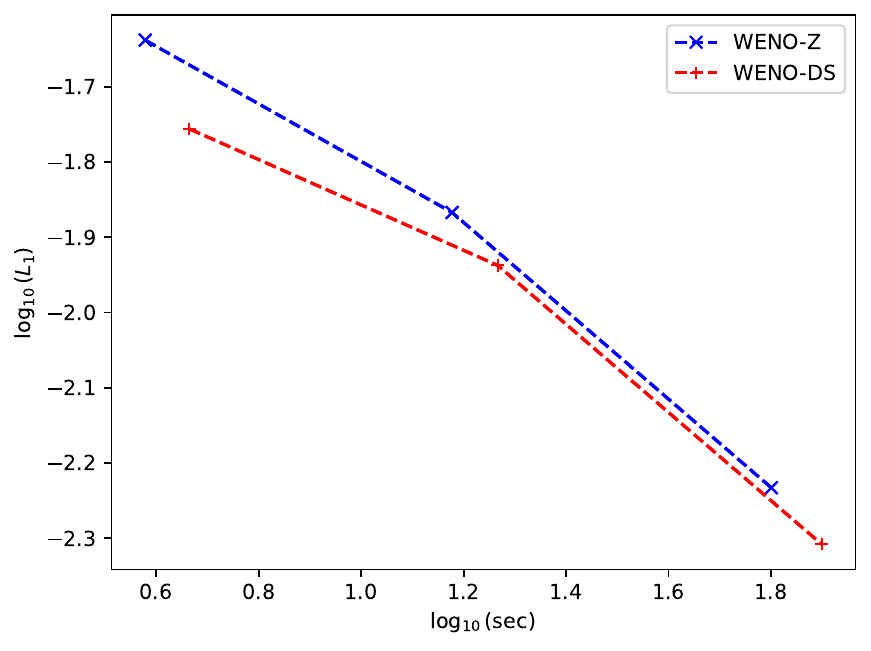}
     \caption{Comparison of computational cost against $L_1$-error of the solution of the Riemann problem with the initial condition \eqref{eq:Riemann_IC_2}.}
\label{fig:Riemann_2_cost}
\end{figure}

Finally, we want to compare the computational cost of WENO-DS compared to the original WENO scheme in solving the problem shown in Figure~\ref{fig:Riemann_2_cost}.
Using a logarithmic scale, we plot the computation time against the $L_1$ error averaged over the four variables $\rho$, $u$, $v$, $p$.    

It should be noted that if we were to test the method on another unseen test problem using the initial data from the previously described range, we would obtain very similar error improvements in those cases.

\subsection{Configuration 3} \label{sec:S5.2}

This is the configuration with four shock waves:
$\overleftarrow{S}_{21}$, $\overleftarrow{S}_{32}$, $\overleftarrow{S}_{34}$, $\overleftarrow{S}_{41}$.
According to \cite{schulz1993}, in this case we have the following equations that must be satisfied:
\begin{equation}  \label{eq:relations_R3}
\begin{split}
    &u_2-u_1 = \Psi_{21}, \quad u_3-u_4 = \Psi_{34}, \quad u_3=u_2, \quad u_4=u_1, \\
    &v_4-v_1 = \Psi_{41}, \quad v_3-v_2 = \Psi_{32}, \quad v_2=v_1, \quad v_3=v_4
\end{split}
\end{equation}
and for polytropic gas the equations
\begin{equation} \label{eq:relations_R3_poly}
    \rho_l/\rho_r = \Pi_{lr} \quad \text{for} \quad (l,r) \in \{(2,1),(3,4),(3,2),(4,1)\}
\end{equation}
are added.
This gives the compatibility conditions $\Psi_{21}=\Psi_{34}$ and $\Psi_{41}=\Psi_{32}$.
Furthermore, we have $\rho_2=\rho_4$, $p_2=p_4$ and $u_2-u_1=v_4-v_1$. 

In this case, we use them for creating the data set values
\begin{equation} \label{eq:par_range_R3}
\begin{split}
      \rho_1 \in \mathcal{U}[1, 2], \quad \rho_2 \in & \mathcal{U}[0.5, 1], \quad p_1 \in \mathcal{U}[1, 2], \\
  \quad u_1 \in \mathcal{U}[-0.25,0.25],  \quad &v_1 = u_1, \quad \gamma \in (1.1, 1.67)
\end{split}
\end{equation}
and for the other values we use the relations \eqref{eq:relations_R3}, \eqref{eq:relations_R3_poly} with \eqref{eq:relations_phi_psi} and \eqref{eq:relations_pi}. 
Similar to the previous example, we compute the reference solutions using the WENO-Z method on a grid $I\times J = 400\times400$ space points up to the final time $T\in\mathcal{U}[0.1,0.3]$ and create the data set consisting of $50$ reference solutions.

\begin{figure}[htb] 
\centering
    \includegraphics[width=0.7\textwidth]{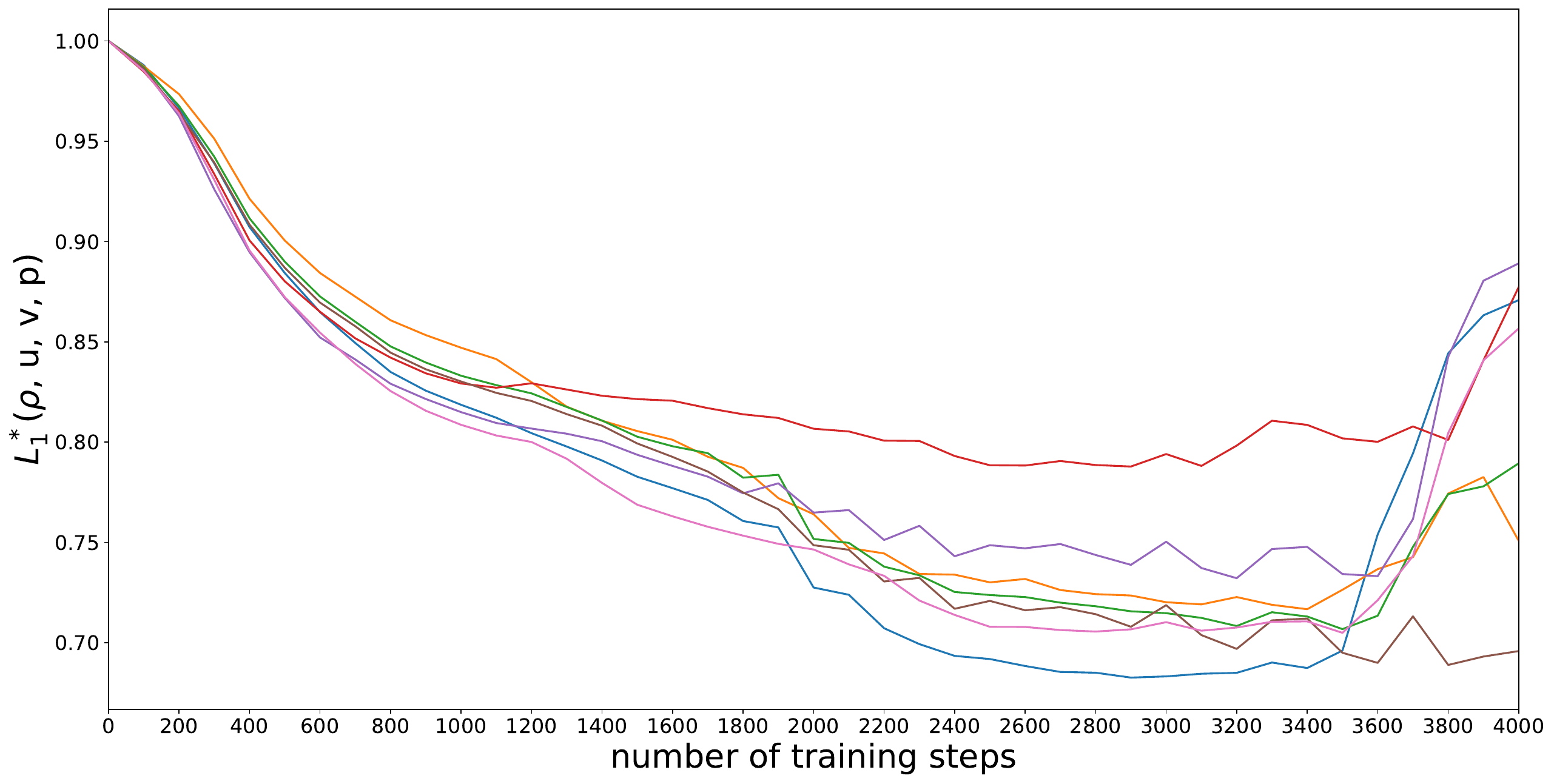}
     \caption{The values \eqref{eq:Euler_validation_adjusted} for different validation problems evaluated each $100$ training steps.}
\label{fig:validation_Riemann_3}
\end{figure}

We proceed with training as described in the previous section, using the same neural network structure as shown in Figure~\ref{fig:structure_1}. 
Again, we train only on the discretization $I \times J = 100 \times 100$ space steps. 
We run the training for $4000$ training steps and plot the evolution of the validation metrics \eqref{eq:Euler_validation_adjusted} for the validation problems in Figure~\ref{fig:validation_Riemann_3}.
We measured the error every $100$ training steps and at the randomly chosen final time $T$.  
Based on this, we choose the final model from training step $3200$ and present the results for the test problem with $\gamma = 1.4$, $T=0.3$, and initial condition 
\cite{kurganov2002}
\begin{equation}  \label{eq:Riemann_IC_3}
  (\rho,u, v,p) = \begin{cases} 
    (1.5, 0, 0, 1.5)  \quad &x > 0.5 \quad  \text{and} \quad  y > 0.5,\\
     (0.5323. 1.206, 0, 0.3) \quad &x < 0.5 \quad  \text{and} \quad  y > 0.5,\\
      (0.138, 1.206, 1.206, 0.029) \quad &x < 0.5 \quad  \text{and} \quad  y < 0.5,\\
       (0.5323, 0, 1.206, 0.3) \quad &x > 0.5 \quad  \text{and} \quad  y < 0.5.\\
\end{cases} 
\end{equation}

\begin{figure}[h]
\begin{subfigure}{.32\textwidth}
\centering
\includegraphics[width=\linewidth]{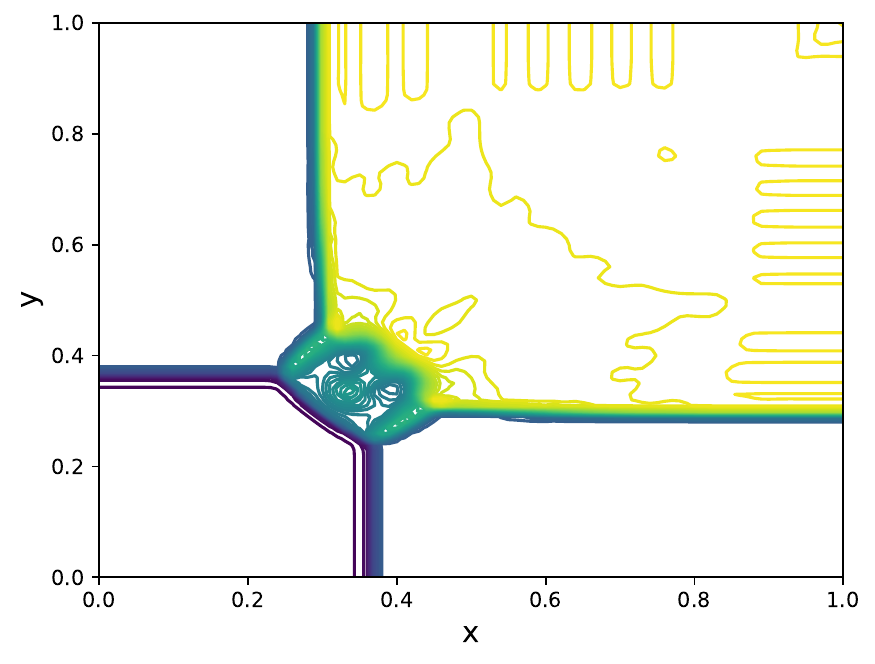}
\caption{WENO-DS}\label{fig:R3_WENODS}
\end{subfigure}
\begin{subfigure}{.32\textwidth}
\centering
\includegraphics[width=\linewidth]{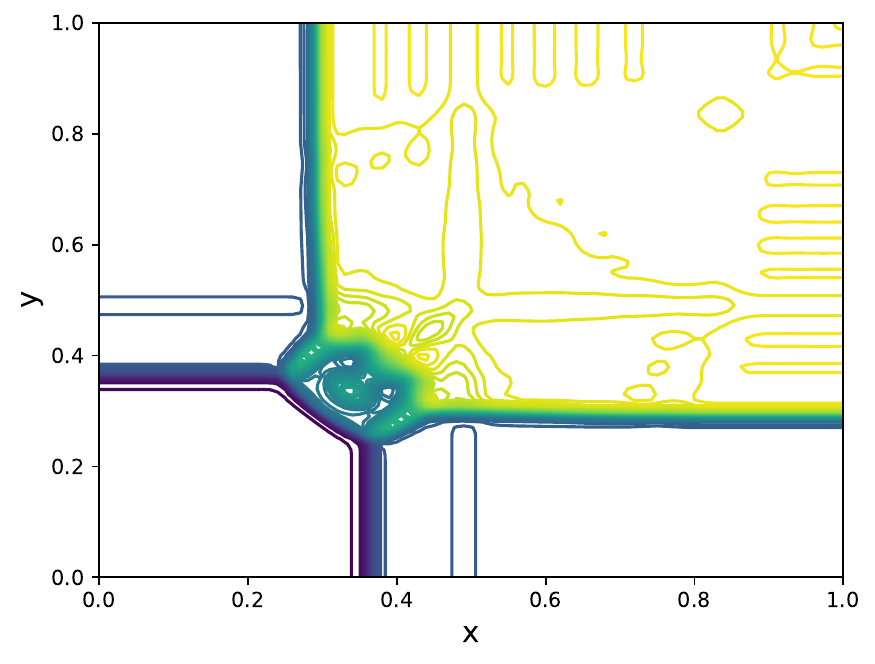} 
\caption{WENO-Z}\label{fig:R3_WENOZ}
\end{subfigure}
\begin{subfigure}{.32\textwidth}
\centering
\includegraphics[width=\linewidth]{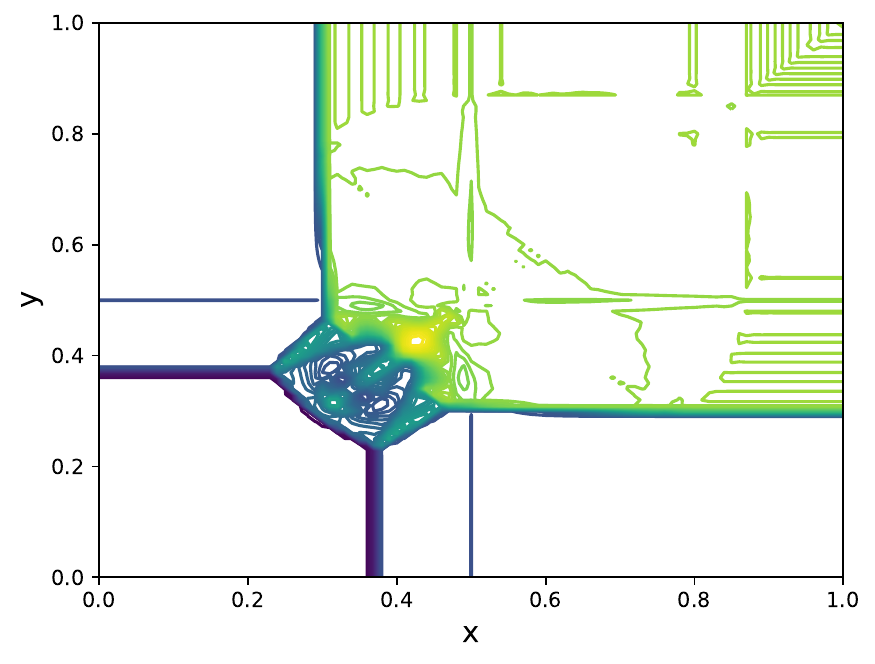} 
\caption{reference solution}\label{fig:R3_reference}
\end{subfigure}
\caption{Density contour plot for the solution of the Riemann problem with the initial condition \eqref{eq:Riemann_IC_3}, $I \times J = 100 \times 100$, $T=0.3$.} \label{fig:Riemann_3}
\end{figure}

\begin{figure}[h!]
\centering
\begin{subfigure}{.4\textwidth}
\centering
\includegraphics[width=\linewidth]{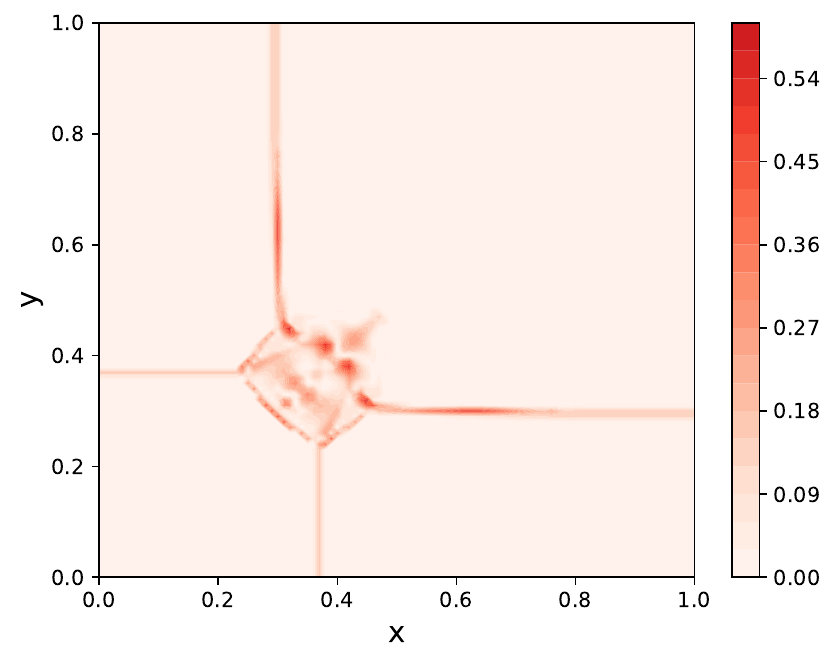}
\caption{WENO-DS}\label{fig:R3_WENODS_error}
\end{subfigure}
\begin{subfigure}{.4\textwidth}
\centering
\includegraphics[width=\linewidth]{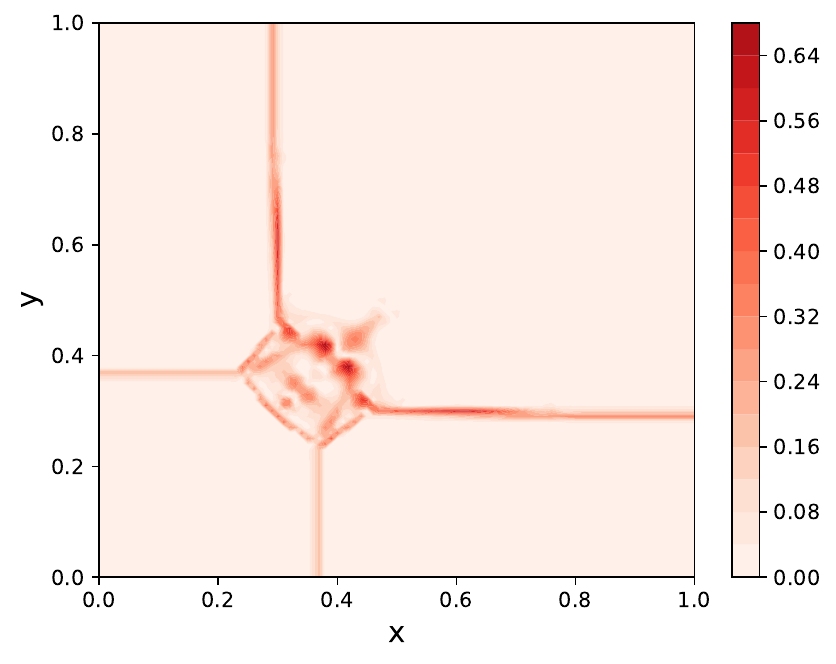} 
\caption{WENO-Z}\label{fig:R3_WENOZ_error}
\end{subfigure}
\caption{Absolute pointwise errors for the density solution of the Riemann problem with the initial condition \eqref{eq:Riemann_IC_3}, $I \times J = 100 \times 100$, $T=0.3$.} \label{fig:Riemann_3_error}
\end{figure}

We compare the results in Table~\ref{tab:Riemann_3}. 
As can be seen, we achieve a large error improvement for all discretizations listed. 
The density contour plots can be found in Figure~\ref{fig:Riemann_3} and the absolute pointwise errors for the density in Figure~\ref{fig:Riemann_3_error}. 
Here it can be seen that the error of WENO-DS is significantly lower in the areas of the shock contacts.

\begin{table}[h] 
 \centering
    \scalebox{0.7}{ \begin{tabular}{|c|c|c|c|c|c|c|c|c|c|}
    \hline
    \multicolumn{1}{|c|}{$I \times J$}&\multicolumn{3}{|c|}{\ $50 \times 50$}&\multicolumn{3}{|c|}{\ $100 \times 100$} &\multicolumn{3}{|c|}{\ $200 \times 200$}\\
    \hline
        & WENO-Z &  WENO-DS & ratio & WENO-Z &  WENO-DS & ratio& WENO-Z &  WENO-DS & ratio\\
    \hline
    \hline
 $\rho$ &   0.038682 &  0.027906 &  1.39  &  0.019232 &  0.012817 &  1.50  &  0.007454 &  0.004657 &  1.60 \\    \hline
$u$ &   0.034692 &  0.027638 &  1.26   &  0.019588 &  0.015043 &  1.30  &  0.008249 &  0.005810 &  1.42 \\     \hline
$v$ &    0.034692 &  0.027638 &  1.26   &  0.019588 &  0.015043 &  1.30  &  0.008249 &  0.005810 &  1.42 \\     \hline
 $p$ &     0.038920 &  0.030888 &  1.26   &  0.018666 &  0.014041 &  1.33  &  0.007275 &  0.005001 &  1.45 \\    \hline
\end{tabular} }
    \caption{Comparison of $L_1$ error of WENO-Z and WENO-DS methods for the solution of the Euler system with the initial condition \eqref{eq:Riemann_IC_3} for different spatial discretizations, $T=0.3$.}
\label{tab:Riemann_3}
\end{table}

We also compare the weights $\omega^Z_m$, $m = 0,1,2$ \eqref{eq:WENOZ} and the updated weights $\omega^{DS}_m$, $m = 0,1,2$ with the improved smoothness indicators \eqref{eq:betas_DS_cons}.
We plot these weights, corresponding to the positive part of a flux $\hat{f}^+$ from the flux splitting, using WENO-Z and WENO-DS for the previous test problem at the final time $T=0.3$. 
Since we apply the principle dimension-by-dimension, we present the weights only for the approximation of the flux $F(U)$. 
For the approximations of the flux $G(U)$, we could obtain these weights in this example using symmetry.
As can be seen, WENO-DS is much better at localizing the shock from the other direction as well, which has a significant impact on error improvement.

\begin{figure}[h] 
\begin{subfigure}{0.32\textwidth}
    \includegraphics[width=\textwidth]{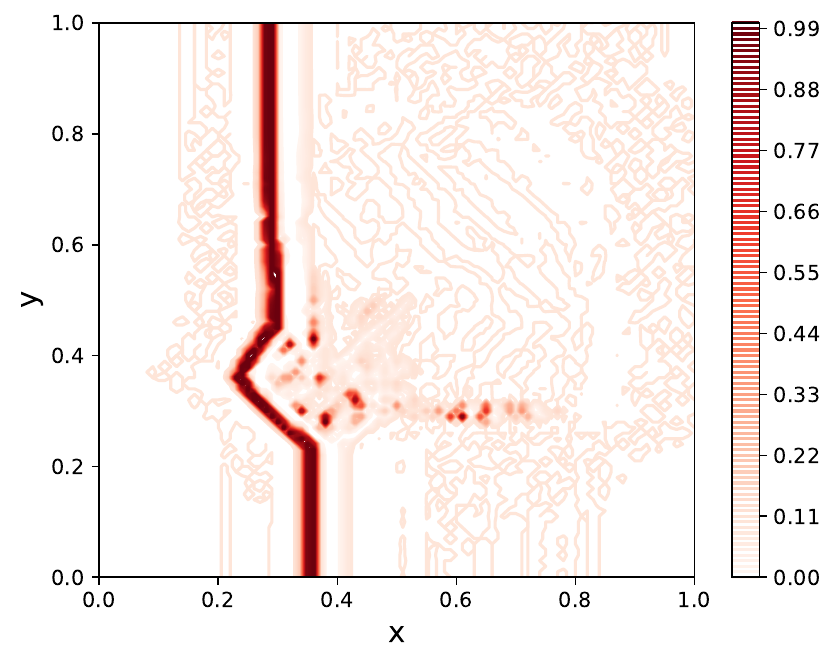}
    \caption{WENO-DS, $\omega^{DS}_0$}
    \label{fig:R3_WENODS_omega_0}
\end{subfigure}
\begin{subfigure}{0.32\textwidth}
     \includegraphics[width=\textwidth]{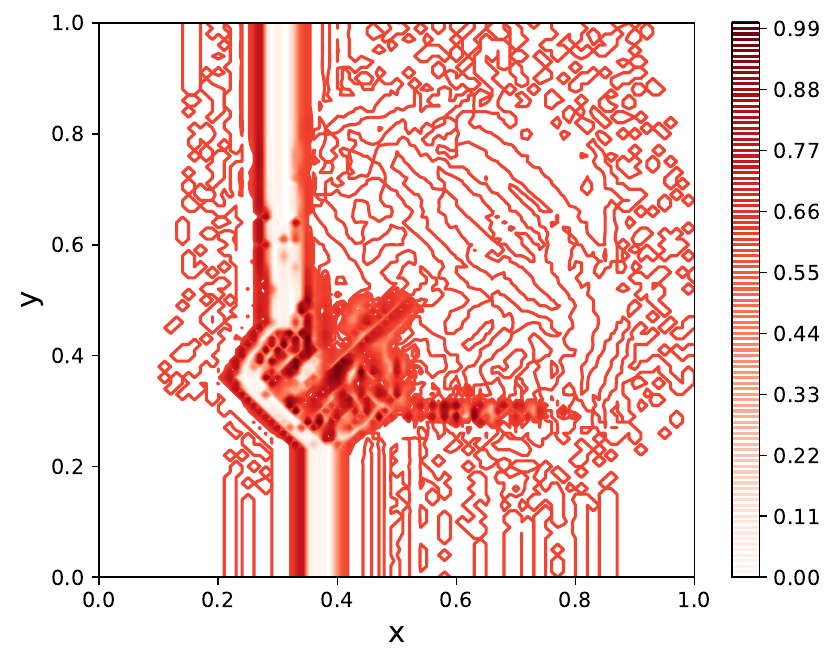}
    \caption{WENO-DS, $\omega^{DS}_1$}
        \label{fig:R3_WENODS_omega_1}
\end{subfigure}
\begin{subfigure}{0.32\textwidth}
     \includegraphics[width=\textwidth]{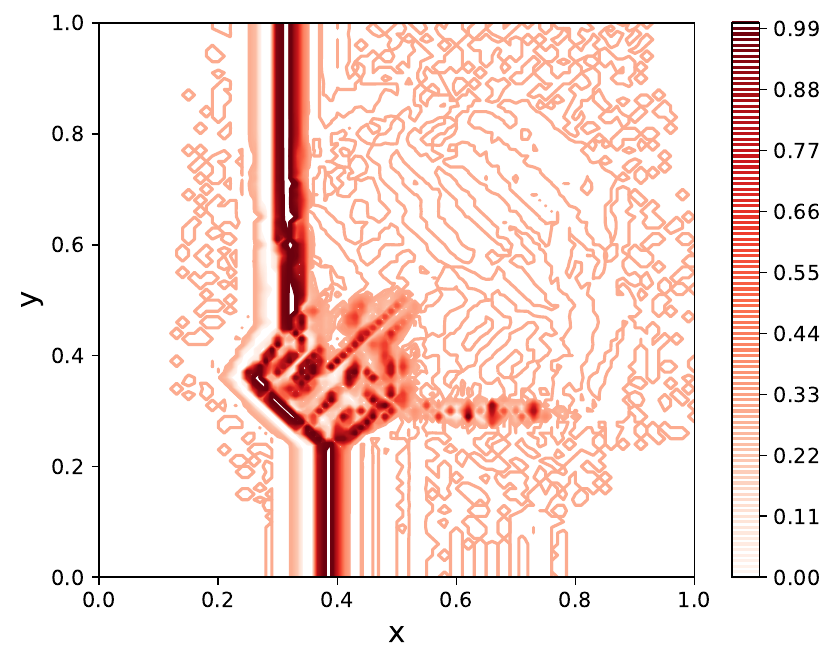}
    \caption{WENO-DS, $\omega^{DS}_2$}
        \label{fig:R3_WENODS_omega_2}
\end{subfigure}
\bigskip
\begin{subfigure}{0.32\textwidth}
    \includegraphics[width=\textwidth]{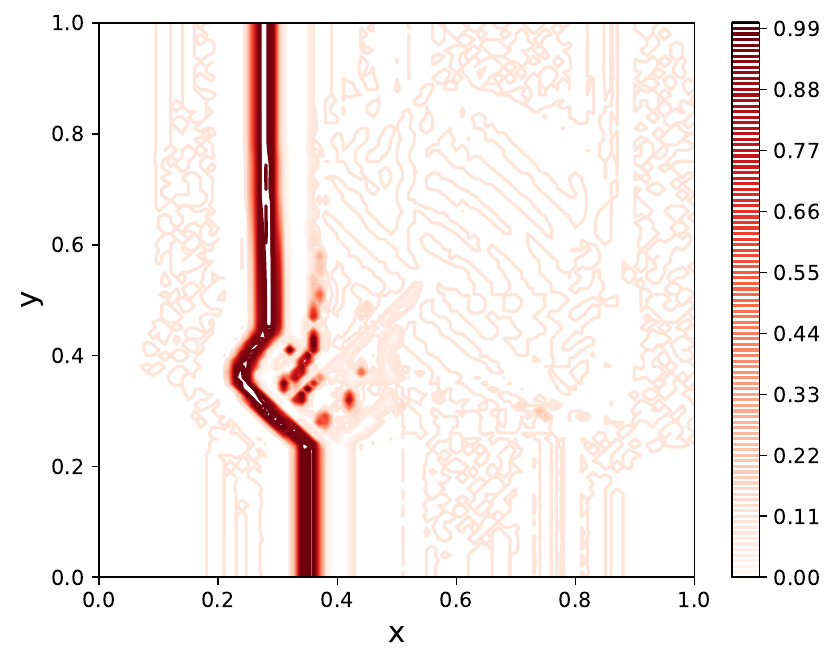}
    \caption{WENO-Z, $\omega^{Z}_0$}
    \label{fig:R3_WENOZ_omega_0}
\end{subfigure}
\begin{subfigure}{0.32\textwidth}
     \includegraphics[width=\textwidth]{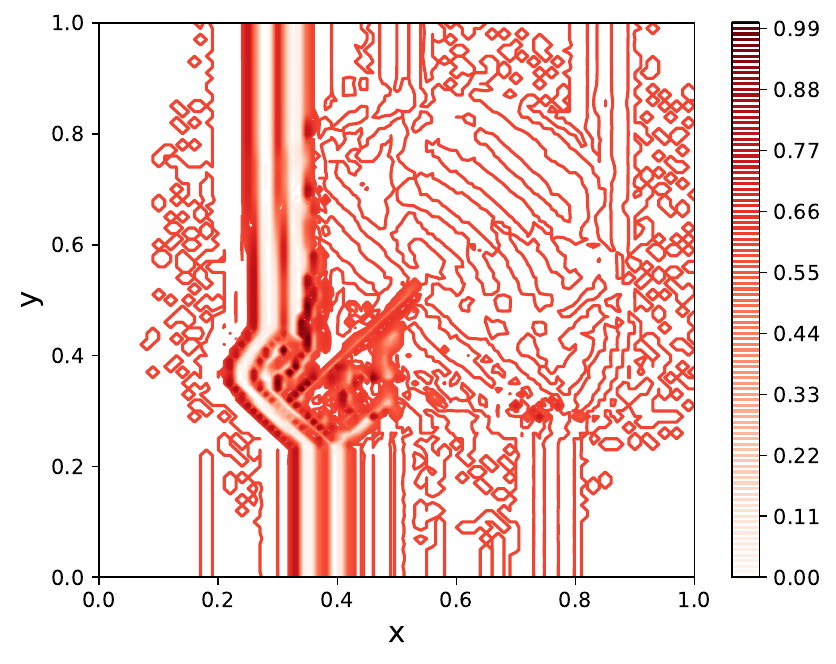}
    \caption{WENO-Z, $\omega^{Z}_1$}\label{fig:R3_WENOZ_omega_1}
\end{subfigure}
\begin{subfigure}{0.32\textwidth}
     \includegraphics[width=\textwidth]{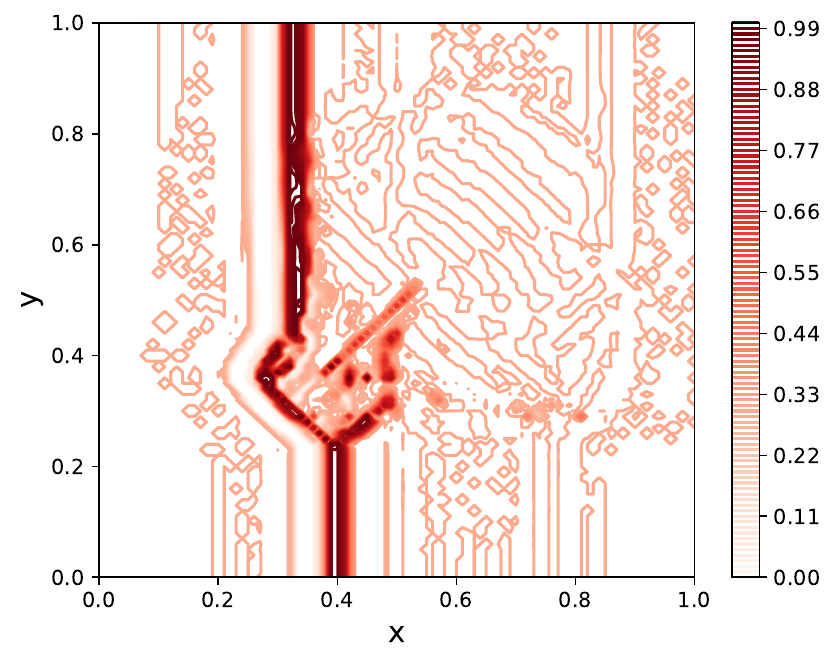}
    \caption{WENO-Z, $\omega^{Z}_2$} \label{fig:R3_WENOZ_omega_2}
\end{subfigure}
\caption{Comparison of the nonlinear weights $\omega^{Z}_m, m = 0,1,2$ and $\omega^{DS}_m$, $m = 0,1,2$.}
\label{fig:R3_omegas}
\end{figure}
\begin{figure}[h!] 
\centering
    \includegraphics[width=0.5\textwidth]{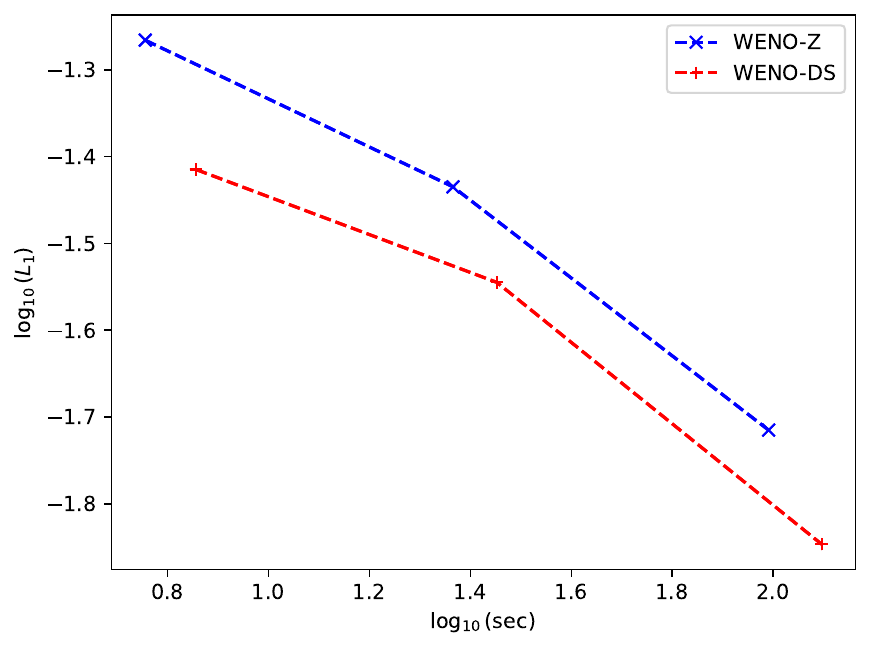}
     \caption{Comparison of computational cost against $L_1$-error of the solution of the Riemann problem with the initial condition \eqref{eq:Riemann_IC_3}.}
\label{fig:Riemann_3_cost}
\end{figure}

Finally, let us compare the computational cost of WENO-DS for the problem shown in Figure~\ref{fig:Riemann_3_cost}. 
We see that WENO-DS is much more computationally intensive compared to WENO-DS.
Again, if we tested the method on the unseen problems, but with the same initial configuration, we would get analogous significant error improvements.

\subsection{Configuration 16} \label{sec:S5.3}

This is the configuration with the combination of rarefaction wave, shock wave and contact discontinuities: $\overleftarrow{R}_{21}$, $J^-_{32}$, $J^+_{34}$, $\overrightarrow{S}_{41}$.
As shown in \cite{schulz1993}, the following relations must hold for this case
\begin{equation}  \label{eq:relations_R16}
\begin{split}
    u_1-u_2 = &\Phi_{21}, \quad u_3=u_4=u_1, \\
    v_4-v_1 = \Psi_{41}, \quad &v_3=v_2 = v_1, \quad p_1 < p_2 = p_3 = p_4
\end{split}
\end{equation}
and for polytropic gas we add the equation \eqref{eq:relations_R2_poly} for a rarefaction and \eqref{eq:relations_R3_poly} for a shock wave between the $l$th and $r$th quadrants.

For our data set we use the values
\begin{equation} \label{eq:par_range_R16}
\begin{split}
      \rho_4 \in \mathcal{U}[1, 2], \quad \rho_3 \in \mathcal{U}[0.5, \rho_4], &\quad p_1 \in \mathcal{U}[0.3,1], \quad p_2 \in \mathcal{U}[1,1.5], \\
  \quad u_1 \in \mathcal{U}[-0.25,0.25],  \quad &v_1 = u_1, \quad \gamma \in (1.1, 1.67)
\end{split}
\end{equation}
To compute the data set consisting of $50$ reference solutions, we use the WENO-Z method on a grid $I \times J = 400 \times400$ space points up to the final time $T\in\mathcal{U}[0.1,0.2]$.

\begin{figure}[htb] 
\centering
    \includegraphics[width=0.7\textwidth]{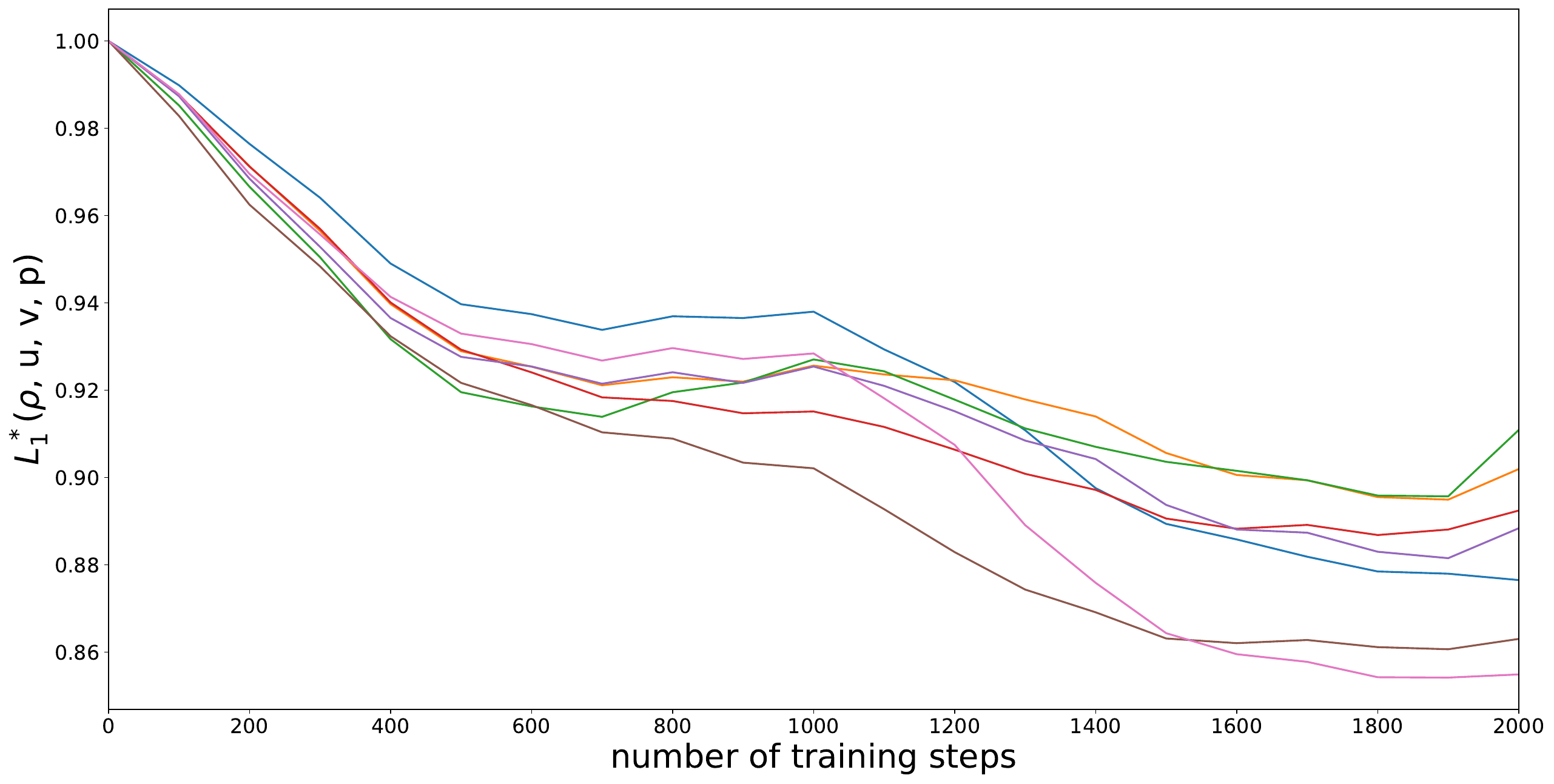}
     \caption{The values \eqref{eq:Euler_validation_adjusted} for different validation problems evaluated each $100$ training steps.}
\label{fig:validation_Riemann_16}
\end{figure}

We train the CNN with the structure shown in Figure~\ref{fig:structure_1} as in the previous examples on the discretization $I \times J = 100 \times 100$ space steps for the total number of $2000$ training steps. 
We show the evolution of the validation metrics \eqref{eq:Euler_validation_adjusted} in Figure~\ref{fig:validation_Riemann_16} and choose the model from training step $1900$.

We test the trained WENO-DS on a test problem \cite{kurganov2002} with $\gamma = 1.4$, $T=0.2$ and the initial condition
\begin{equation}  \label{eq:Riemann_IC_16}
  (\rho,u, v,p) = \begin{cases} 
    (0.5313, 0.1, 0.1, 0.4)  \quad &x > 0.5 \quad  \text{and} \quad  y > 0.5,\\
     (1.0222, -0.6179, 0.1, 1) \quad &x < 0.5 \quad  \text{and} \quad  y > 0.5,\\
      (0.8, 0.1, 0.1, 1) \quad &x < 0.5 \quad  \text{and} \quad  y < 0.5,\\
       (1, 0.1, 0.8276, 1) \quad &x > 0.5 \quad  \text{and} \quad  y < 0.5.\\
\end{cases} 
\end{equation}

\begin{figure}[h!]
\begin{subfigure}{.32\textwidth}
\centering
\includegraphics[width=\linewidth]{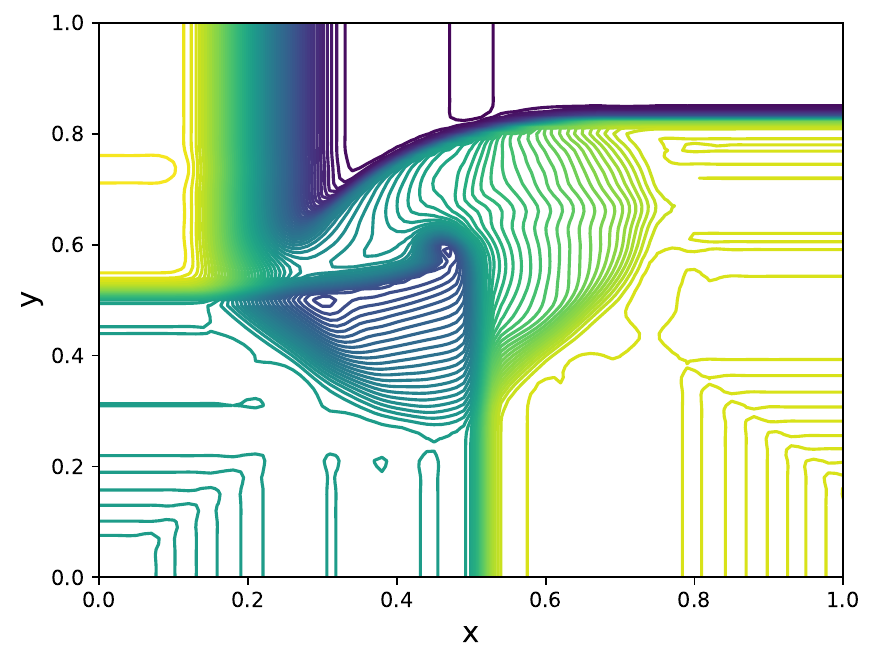}
\caption{WENO-DS}
\label{fig:R16_WENODS}
\end{subfigure}
\begin{subfigure}{.32\textwidth}
\centering
\includegraphics[width=\linewidth]{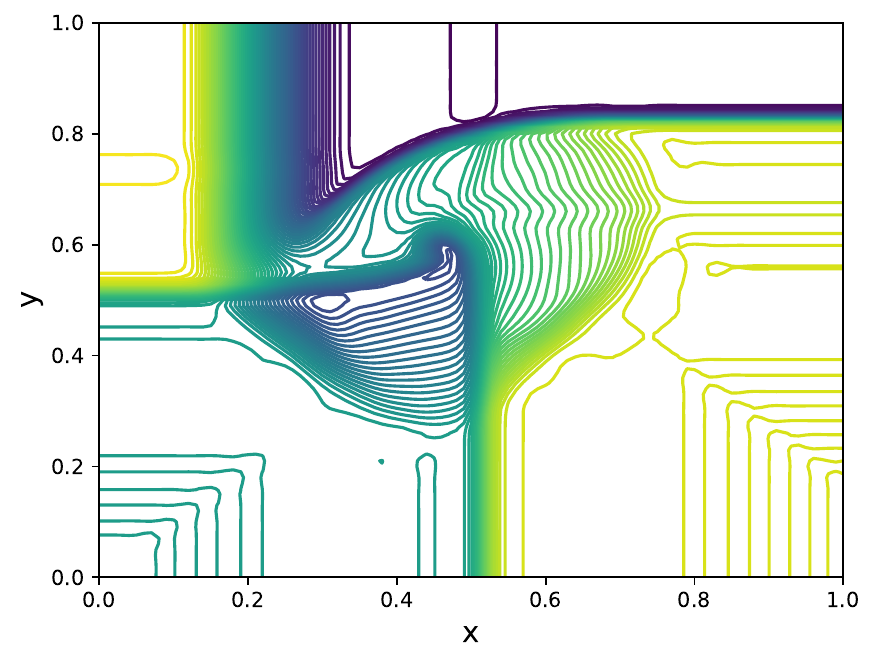} 
\caption{WENO-Z}
\label{fig:R16_WENOZ}
\end{subfigure}
\begin{subfigure}{.32\textwidth}
\centering
\includegraphics[width=\linewidth]{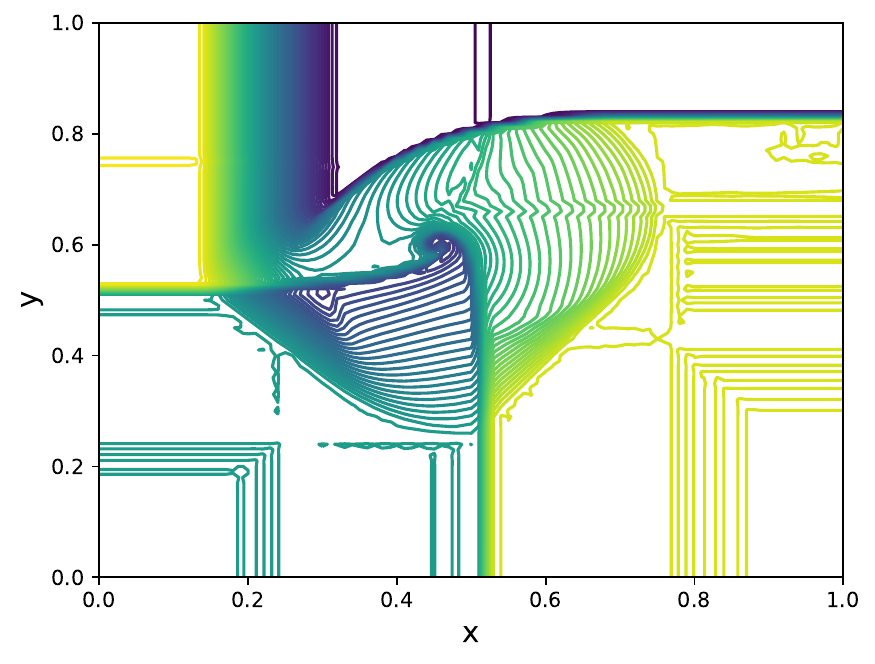} 
\caption{reference solution}
\label{fig:R16_reference}
\end{subfigure}
\caption{Density contour plot for the solution of the Riemann problem with the initial condition \eqref{eq:Riemann_IC_16}, $I \times J = 100 \times 100$, $T=0.2$.} \label{fig:Riemann_16}
\end{figure}

\begin{figure}[h!]
\centering
\begin{subfigure}{.4\textwidth}
\centering
\includegraphics[width=\linewidth]{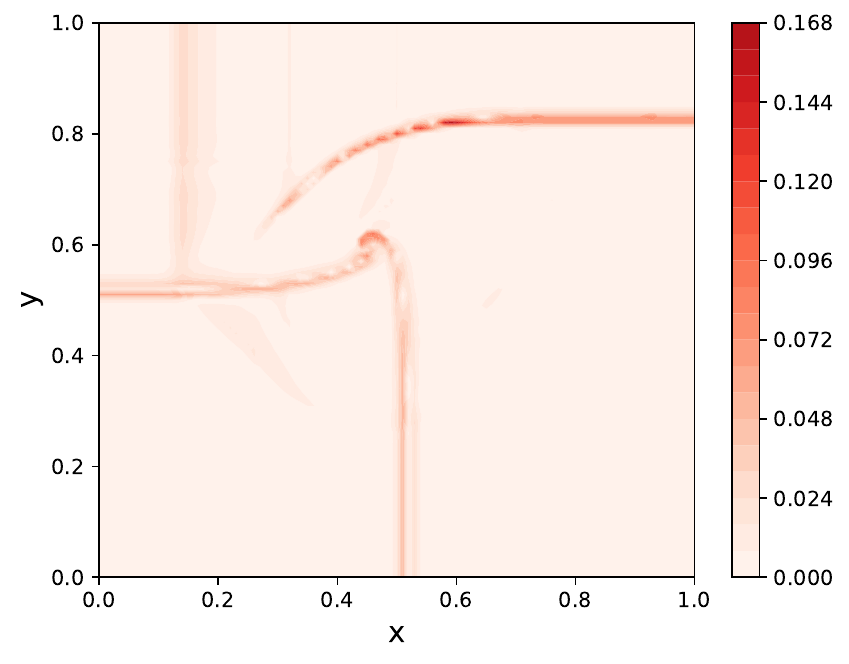}
\caption{WENO-DS}
\label{fig:R16_WENODS_error}
\end{subfigure}
\begin{subfigure}{.4\textwidth}
\centering
\includegraphics[width=\linewidth]{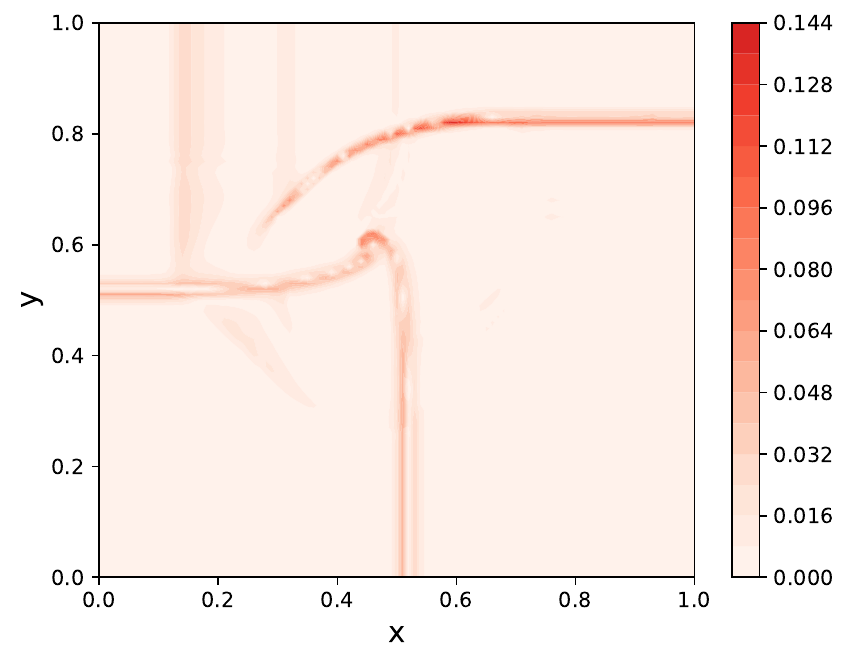} 
\caption{WENO-Z}
\label{fig:R16_WENOZ_error}
\end{subfigure}
\caption{Absolute pointwise errors for density solution of the Riemann problem with the initial condition \eqref{eq:Riemann_IC_16}, $I \times J = 100 \times 100$, $T=0.2$.} \label{fig:Riemann_16_error}
\end{figure}

We compare the results in Table~\ref{tab:Riemann_16} and the density contour plots can be found in Figure~\ref{fig:Riemann_16}. 
As can be seen, WENO-DS outperforms WENO-Z and has smaller $L_1$ errors in all cases.
In addition, we plot the absolute pointwise errors for the density solution and show them in Figure~\ref{fig:Riemann_16_error}. 

For another unseen test problem with the same initial configurations, we would again obtain analogous significant error improvements.

\begin{table}[h] 
 \centering
    \scalebox{0.7}{ \begin{tabular}{|c|c|c|c|c|c|c|c|c|c|}
    \hline
    \multicolumn{1}{|c|}{$I \times J$}&\multicolumn{3}{|c|}{\ $50 \times 50$}&\multicolumn{3}{|c|}{\ $100 \times 100$} &\multicolumn{3}{|c|}{\ $200 \times 200$}\\
    \hline
        & WENO-Z &  WENO-DS & ratio & WENO-Z &  WENO-DS & ratio& WENO-Z &  WENO-DS & ratio\\
    \hline
    \hline
 $\rho$  &  0.010980 &  0.009877 &  1.11  &  0.004834 &  0.004327 &  1.12  &  0.001827 &  0.001624 &  1.12     \\    \hline
$u$   &  0.012464 &  0.011287 &  1.10  &  0.005989 &  0.005326 &  1.12  &  0.002223 &  0.001913 &  1.16    \\     \hline
$v$   &  0.015020 &  0.013932 &  1.08  &  0.006609 &  0.006172 &  1.07  &  0.002527 &  0.002298 &  1.10   \\     \hline
 $p$   &  0.010594 &  0.009644 &  1.10  &  0.004236 &  0.003820 &  1.11  &  0.001576 &  0.001392 &  1.13      \\    \hline
\end{tabular} }
    \caption{Comparison of $L_1$ error of WENO-Z and WENO-DS methods for the solution of the Euler system with the initial condition \eqref{eq:Riemann_IC_16} for different spatial discretizations, $T=0.2$.}
\label{tab:Riemann_16}
\end{table}

\subsection{Configuration 11 and Configuration 19} \label{sec:S5.4}

In the previous sections, we trained three WENO-DS methods for three different types of configurations. 
We denote by WENO-DS (C2), WENO-DS (C3), and WENO-DS (C16) the methods from Sections~\ref{sec:S5.1}, \ref{sec:S5.2}, and \ref{sec:S5.3}, respectively.
In this section, we test these methods on the unseen problems containing the combination of rarefaction wave, shock wave, and contact discontinuities. 
First, we consider Configuration~11 ($\overleftarrow{S}_{21}$, $J^+_{32}$, $J^+_{34}$, $\overleftarrow{S}_{41}$) with the test problem with $\gamma = 1.4$, $T=0.3$, and the initial condition \cite{kurganov2002}
\begin{equation}  \label{eq:Riemann_IC_11}
  (\rho,u, v,p) = \begin{cases} 
    (1, 0.1, 0, 1)  \quad &x > 0.5 \quad  \text{and} \quad  y > 0.5,\\
     (0.5313, 0.8276, 0, 0.4) \quad &x < 0.5 \quad  \text{and} \quad  y > 0.5,\\
      (0.8, 0.1, 0, 0.4) \quad &x < 0.5 \quad  \text{and} \quad  y < 0.5,\\
       (0.5313, 0.1, 0.7276, 0.4) \quad &x > 0.5 \quad  \text{and} \quad  y < 0.5.\\
\end{cases} 
\end{equation}
Second, we test the models on the configuration ($J^+_{21}$, $\overleftarrow{S}_{32}$, $J^-_{34}$, $\overrightarrow{R}_{41}$) with the test problem with $\gamma = 1.4$, $T=0.3$ and the initial condition \cite{kurganov2002}
\begin{equation}\label{eq:Riemann_IC_19}
  (\rho,u, v,p) = \begin{cases} 
    (1, 0, 0.3, 1)  \quad &x > 0.5 \quad  \text{and} \quad  y > 0.5,\\
     (2, 0, -0.3, 1) \quad &x < 0.5 \quad  \text{and} \quad  y > 0.5,\\
      (1.0625, 0, 0.2145, 0.4) \quad &x < 0.5 \quad  \text{and} \quad  y < 0.5,\\
       (0.5197, 0, -0.4259, 0.4) \quad &x > 0.5 \quad  \text{and} \quad  y < 0.5.\\
\end{cases} 
\end{equation}
We summarize the results in Tables~\ref{tab:Riemann_11} and \ref{tab:Riemann_19}. 
As can be seen, the method trained on problems containing only rarefaction waves has the worst ability to generalize to unseen problems. 
On the other hand, by using methods trained on problems containing shocks or a combination of contact discontinuities, rarefaction, and shock waves, 
we obtain the error improvements even on unseen problems with different initial configurations.
We would like to emphasize that the test problems in this section are far from the problems included in the training and validation sets. 
This is not only due to the choice of initial data, but also to the combination of rarefaction, shock waves and their direction, and positive and negative contact discontinuities.

\begin{table}[h] 
 \centering
    \scalebox{0.7}{ \begin{tabular}{|c|c||c|c|c|c|c|c|}
    \hline
    \multicolumn{1}{|c|}{}&\multicolumn{7}{|c|}{\ Configuration 11} \\
    \hline
        & WENO-Z &  WENO-DS (C2) & ratio &  WENO-DS (C3)  & ratio &  WENO-DS (C16) & ratio\\
    \hline
    \hline
 $\rho$      &  0.007792 &  0.008000 &  0.97  &   \textbf{0.006783} &  \textbf{1.15}  &    0.007538 &  1.03 \\    \hline
$u$       &  0.008003 &  0.008701 &  0.92  &   0.007846 &  1.02  &    \textbf{0.007840} &  \textbf{1.02}  \\     \hline
$v$       &  0.007692 &  0.008300 &  0.93  &   \textbf{0.007161} &  \textbf{1.07}  &    0.007370 &  1.04  \\     \hline
 $p$       &  0.005883 &  0.006467 &  0.91  &   \textbf{0.005115} &  \textbf{1.15}  &    0.005776 &  1.02  \\    \hline
\end{tabular} }
    \caption{Comparison of $L_1$ error of WENO-Z and WENO-DS methods trained on data in Sections~\ref{sec:S5.1}, \ref{sec:S5.2} and  \ref{sec:S5.3} for the solution of the Euler system with the initial condition \eqref{eq:Riemann_IC_11}, $I \times J = 100 \times 100$, $T=0.3$.}
\label{tab:Riemann_11}
\end{table}

\begin{table}[h] 
 \centering
    \scalebox{0.7}{ \begin{tabular}{|c|c||c|c|c|c|c|c|}
    \hline
    \multicolumn{1}{|c|}{}&\multicolumn{7}{|c|}{\ Configuration 19} \\
    \hline
        & WENO-Z &  WENO-DS (C2) & ratio &  WENO-DS (C3)  & ratio &  WENO-DS (C16) & ratio\\
    \hline
    \hline
 $\rho$ &  0.014844 &  0.014463 &  1.03 &  \textbf{0.013702} &  \textbf{1.08}  &   0.013841 &  1.07  \\    \hline
$u$  &  0.003749 &  \textbf{0.003562} &  \textbf{1.05}  &  0.003689 &  1.02  &   0.003574 &  1.05  \\     \hline
$v$  &  0.009891 &  0.009502 &  1.04  &  0.009791 &  1.01  &   \textbf{0.009245} &  \textbf{1.07} \\     \hline
 $p$  &  0.006123 &  0.005922 &  1.03  &  \textbf{0.005595} &  \textbf{1.09}  &   0.005844 &  1.05  \\    \hline
\end{tabular} }
    \caption{Comparison of $L_1$ error of WENO-Z and WENO-DS methods trained on data in Sections~\ref{sec:S5.1},  \ref{sec:S5.2} and  \ref{sec:S5.3} for the solution of the Euler system with the initial condition \eqref{eq:Riemann_IC_19}, $I \times J = 100 \times 100$, $T=0.3$.}
\label{tab:Riemann_19}
\end{table}



\subsection{Bigger CNN and ability to generalize on unseen configurations} \label{sec:S5.5}

As can be seen from the previous Section~\ref{sec:S5.4}, the models trained using the data from Section~\ref{sec:S5.2} and Section~\ref{sec:S5.3} are able to generalize very well to unseen problems. 
The WENO-DS method is able to properly localize the shocks and discontinuities, leading to a better numerical solution. 
Let us now increase the size of the CNN and use the structures shown in Figures~\ref{fig:structure_2}, 
increasing the size of the receptive field and the number of channels, and Figure~\ref{fig:structure_3}, increasing the number of channels and adding another CNN layer.
Experimentally, we found that only increasing the size of the receptive field and the number of channels leads to similar results as described in the previous sections.
In addition, increasing the receptive field makes the WENO-DS computationally more expensive.
This is because we need to prepare wider inputs for the CNN, which also need to be projected onto the characteristic fields using left eigenvectors, and the matrix multiplications are more expensive here.
On the other hand, if we use the CNN structure described in Figure~\ref{fig:structure_3}, we obtain a trained WENO-DS method that provides a much better numerical solution even for unseen problems with significantly different initial configurations. 

Let us now train the method on two data sets. 
First, we use the dataset from the Section~\ref{sec:S5.2}, train the CNN, and denote the final method as WENO-DS (C3c). 
Second, we train the CNN on the data set from the Section~\ref{sec:S5.3} and denote the final method as WENO-DS (C16c). 
We test the methods on even more different configurations and compare the results in Tables~\ref{tab:Riemann_C3c} and \ref{tab:Riemann_C16c}. 
We use boldface to indicate the configuration on which the method was actually trained. 

With the number of configurations listed in the tables, we cover a wide range of possible combinations of contact discontinuities, rarefaction and shock waves. 
For all of them we use the test examples from the literature, see, e.g.\ \cite{kurganov2002}. 
We treat the possibility with four contact discontinuities with Configuration~6: $J^-_{21}$, $J^-_{32}$, $J^-_{34}$, $J^-_{41}$, 
two contact discontinuities and two rarefaction waves with Configuration~8: $\overleftarrow{R}_{21}$, $J^-_{32}$, $J^-_{34}$, $\overleftarrow{R}_{41}$, two shock waves and two contact discontinuities using Configuration~14: $J^+_{21}$, $\overleftarrow{S}_{32}$, $J^-_{34}$, $\overleftarrow{S}_{41}$ and Configuration~11 from Section~\ref{sec:S5.4}. 
Finally, the combination of contact discontinuities, rarefaction, and shock waves using Configuration~18: $J^+_{21}$, $\overleftarrow{S}_{32}$, $J^+_{34}$, $\overrightarrow{R}_{41}$, and Configuration~19 form the Section~\ref{sec:S5.4}.

\begin{figure}[h!]
\centering
\begin{subfigure}{.32\textwidth}
\centering
\includegraphics[width=\linewidth]{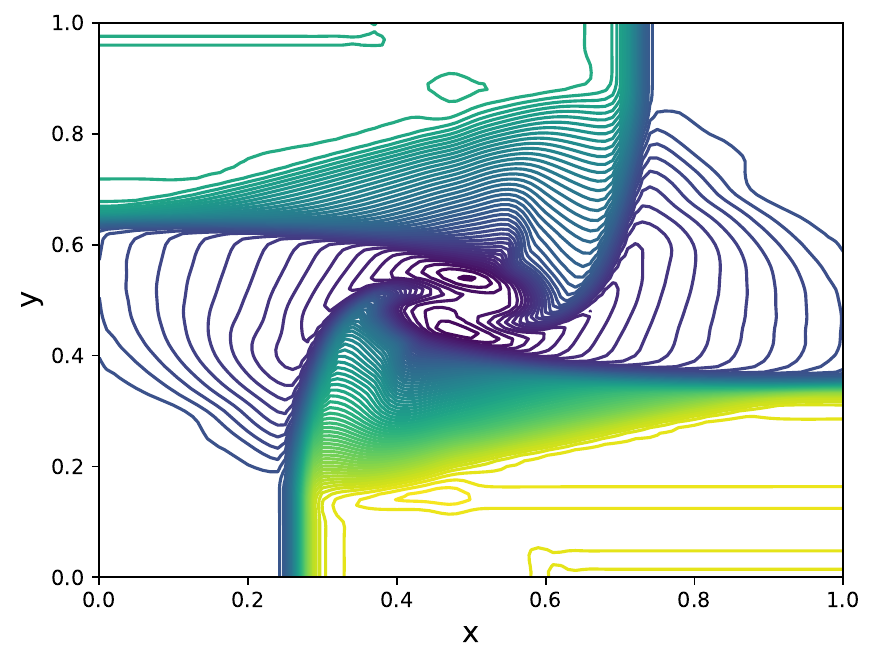}
\caption{WENO-DS}
\label{fig:R6_WENODS}
\end{subfigure}
\begin{subfigure}{.32\textwidth}
\centering
\includegraphics[width=\linewidth]{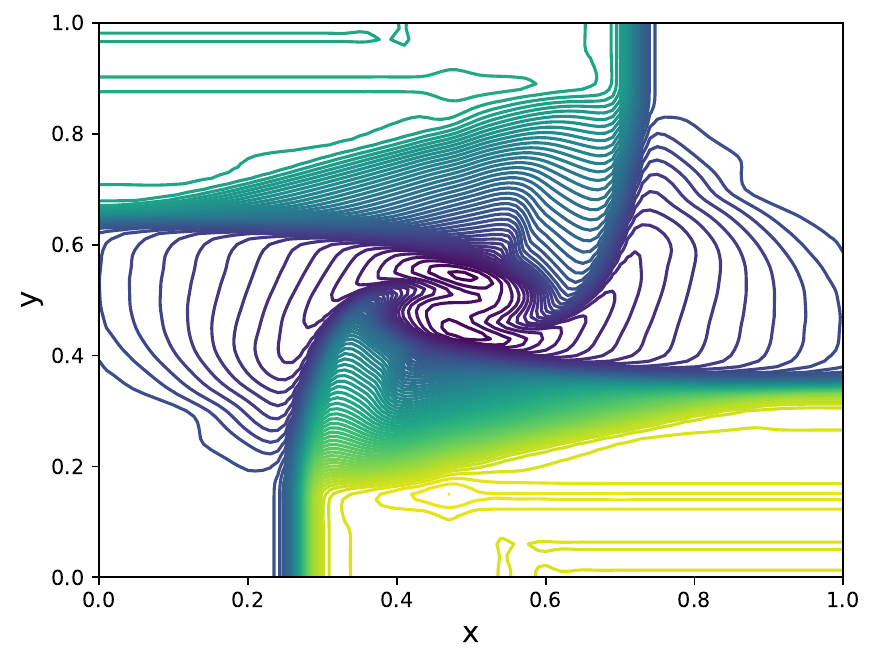} 
\caption{WENO-Z}
\label{fig:R6_WENOZ}
\end{subfigure}
\begin{subfigure}{.32\textwidth}
\centering
\includegraphics[width=\linewidth]{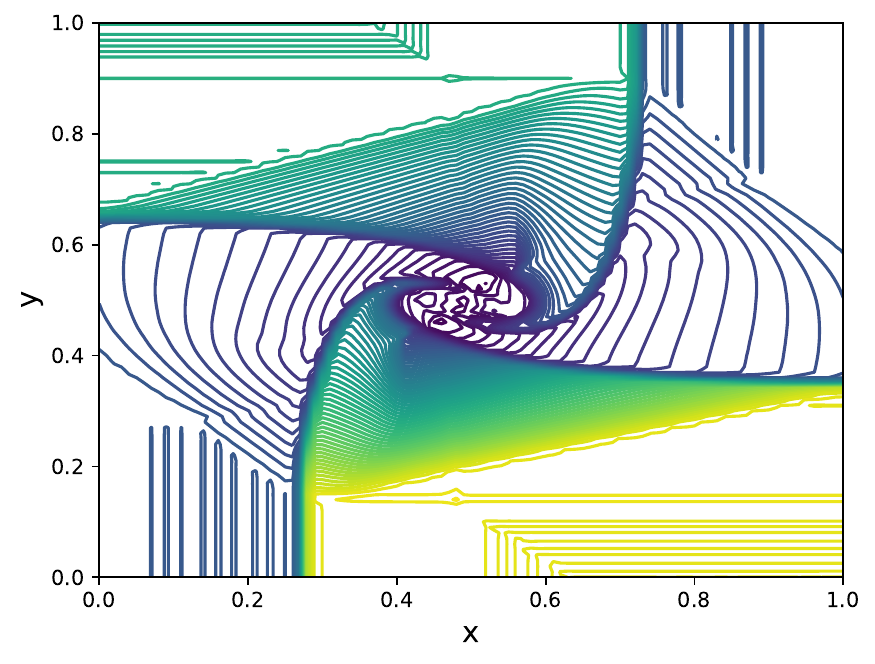} 
\caption{reference solution}
\label{fig:R6_reference}
\end{subfigure}
\begin{subfigure}{.32\textwidth}
\centering
\includegraphics[width=\linewidth]{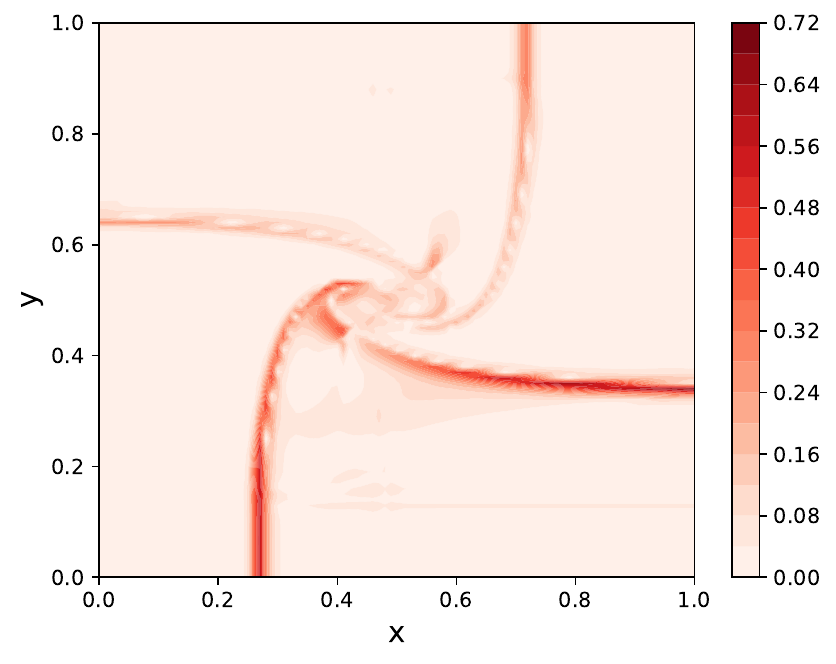} 
\caption{WENO-DS}
\end{subfigure}
\begin{subfigure}{.32\textwidth}
\centering
\includegraphics[width=\linewidth]{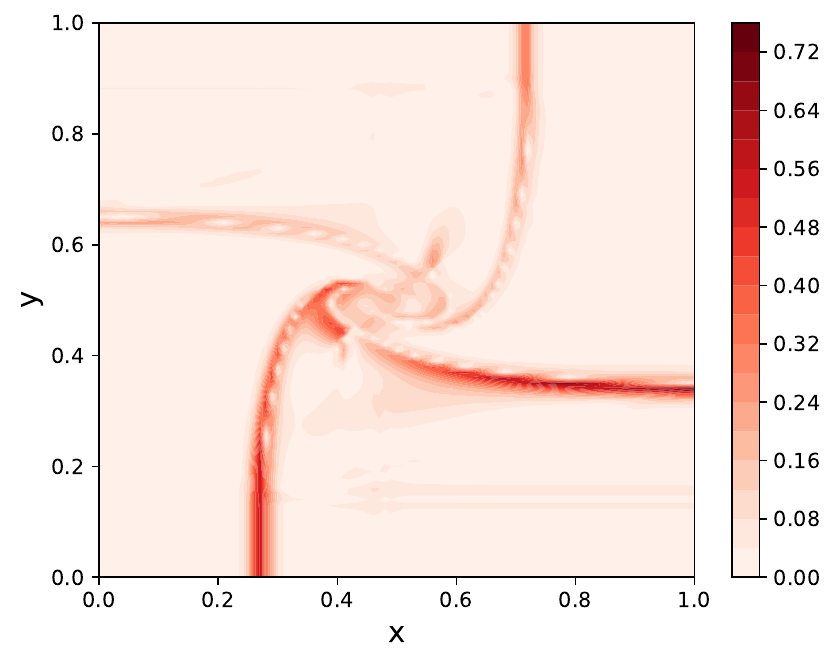} 
\caption{WENO-Z}
\end{subfigure}
\caption{Density contour plots and absolute pointwise errors for the solution of the Riemann problem with initial Configuration~6, $I \times J = 100 \times 100$, $T=0.3$.} \label{fig:Riemann_6}
\end{figure}

\begin{figure}[h!]
\centering
\begin{subfigure}{.32\textwidth}
\centering
\includegraphics[width=\linewidth]{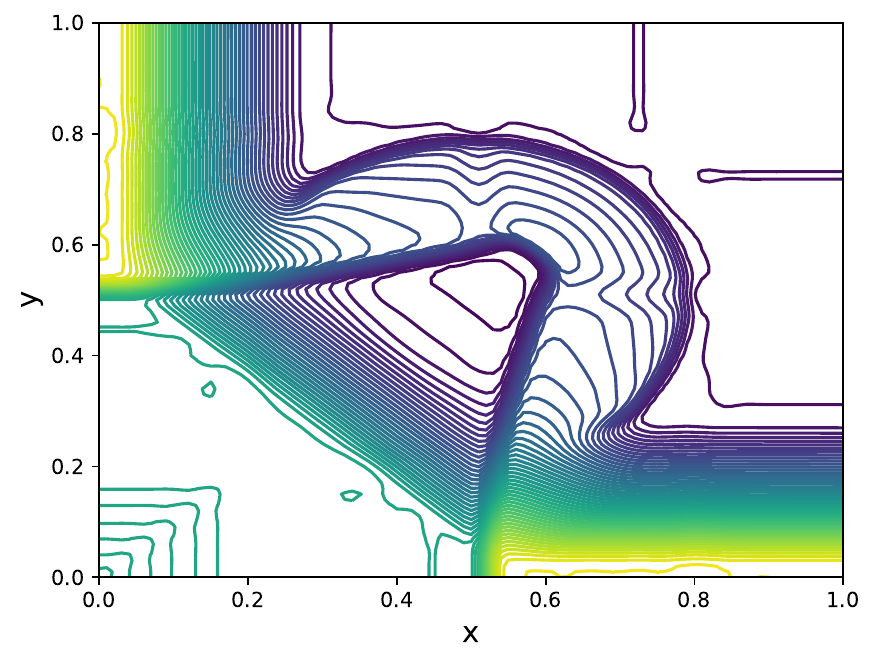}
\caption{WENO-DS}\label{fig:R8_WENODS}
\end{subfigure}
\begin{subfigure}{.32\textwidth}
\centering
\includegraphics[width=\linewidth]{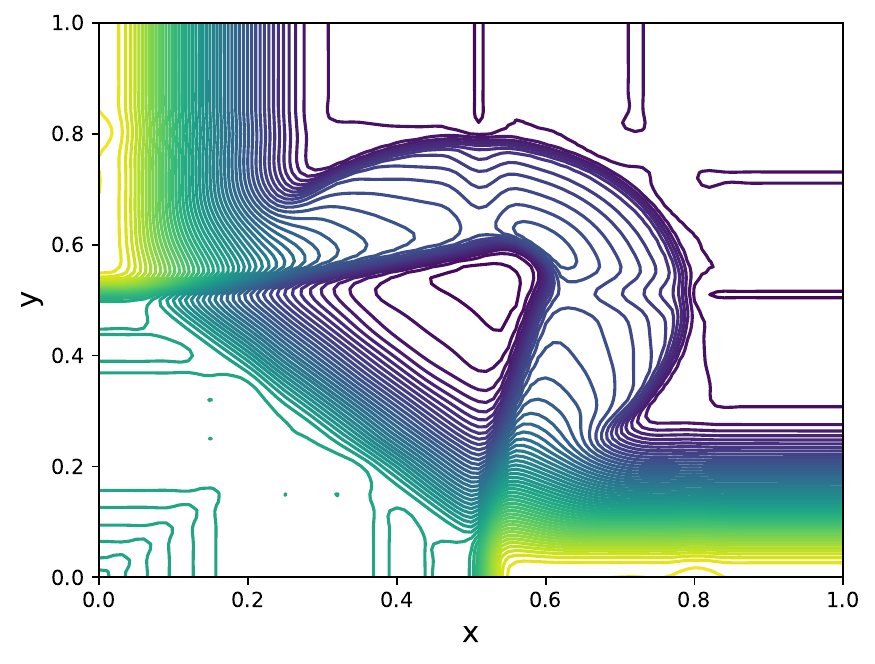} 
\caption{WENO-Z}\label{fig:R8_WENOZ}
\end{subfigure}
\begin{subfigure}{.32\textwidth}
\centering
\includegraphics[width=\linewidth]{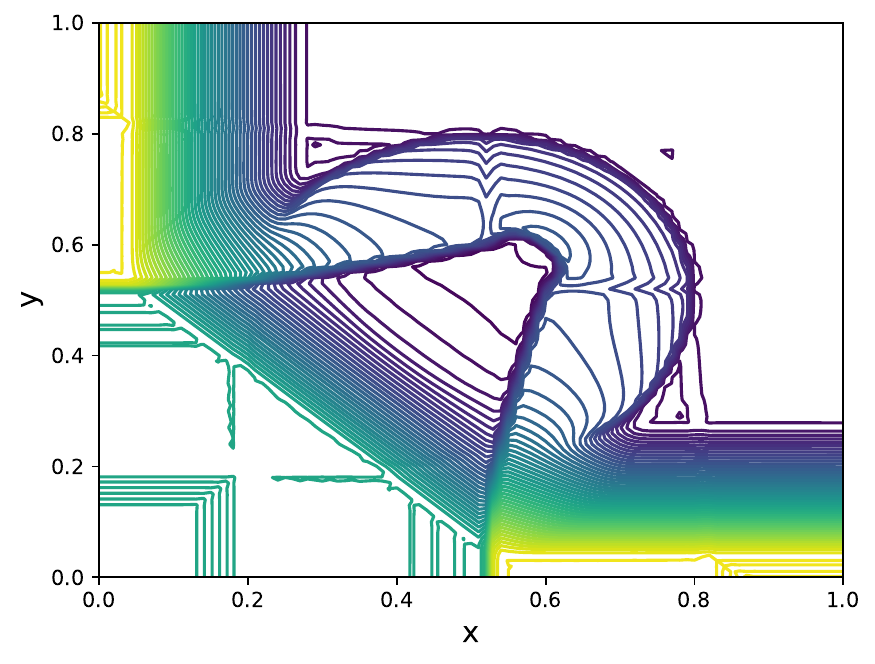} 
\caption{reference solution}\label{fig:R8_reference}
\end{subfigure}
\begin{subfigure}{.32\textwidth}
\centering
\includegraphics[width=\linewidth]{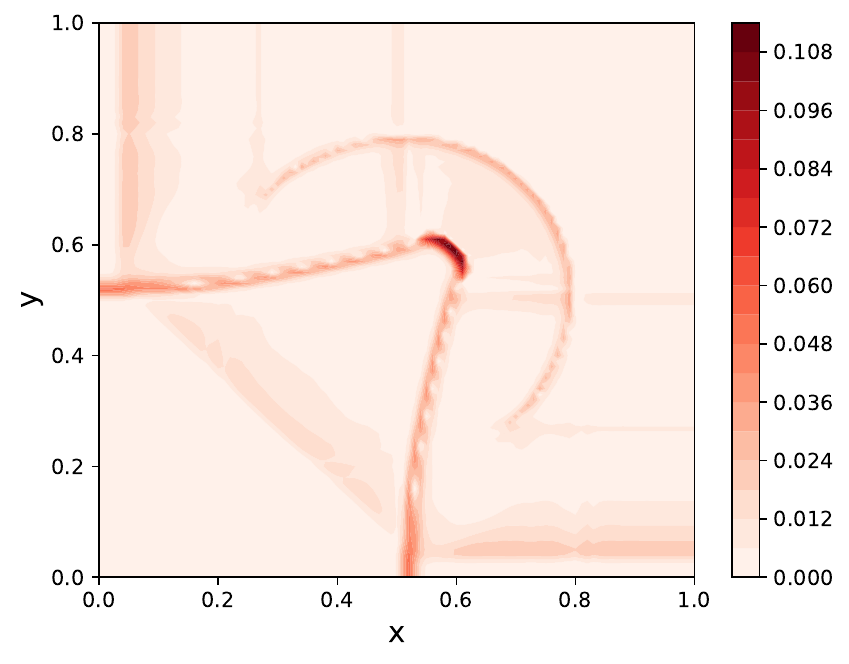} 
\caption{WENO-DS}
\end{subfigure}
\begin{subfigure}{.32\textwidth}
\centering
\includegraphics[width=\linewidth]{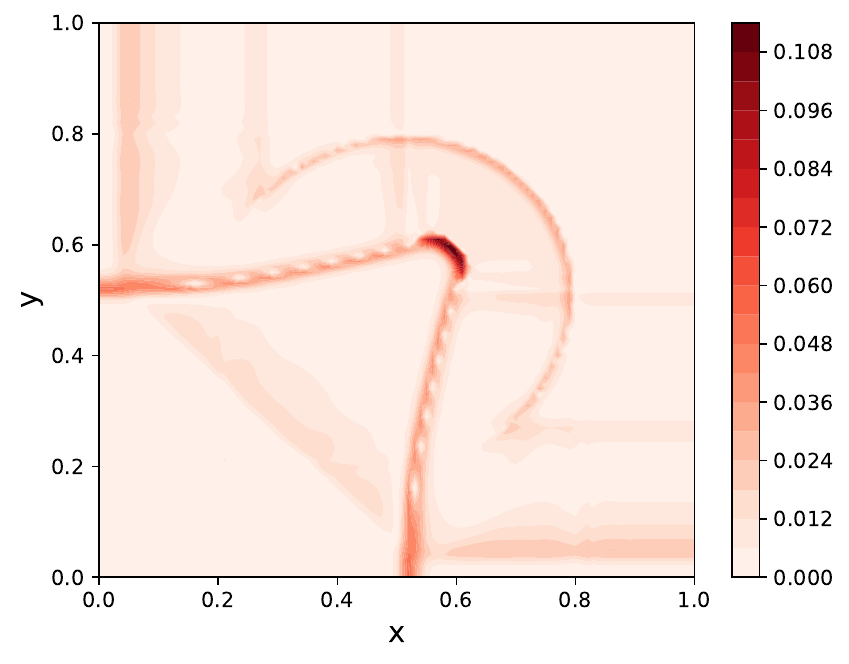} 
\caption{WENO-Z}
\end{subfigure}
\caption{Density contour plots and absolute pointwise errors for the solution of the Riemann problem with initial Configuration~8, $I \times J = 100 \times 100$, $T=0.25$.} \label{fig:Riemann_8}
\end{figure}

\begin{figure}[h!]
\centering
\begin{subfigure}{.32\textwidth}
\centering
\includegraphics[width=\linewidth]{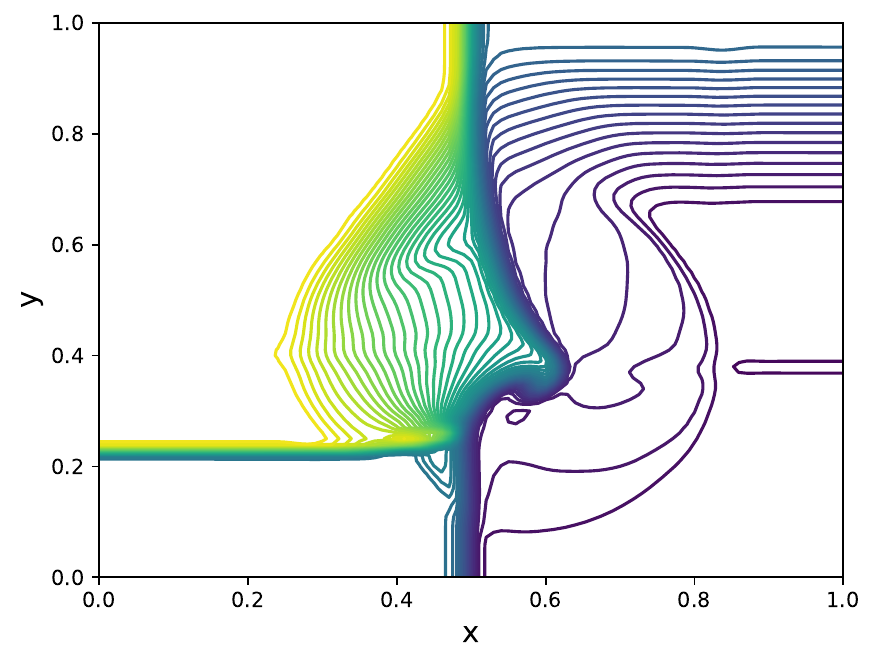}
\caption{WENO-DS}
\label{fig:R19_WENODS}
\end{subfigure}
\begin{subfigure}{.32\textwidth}
\centering
\includegraphics[width=\linewidth]{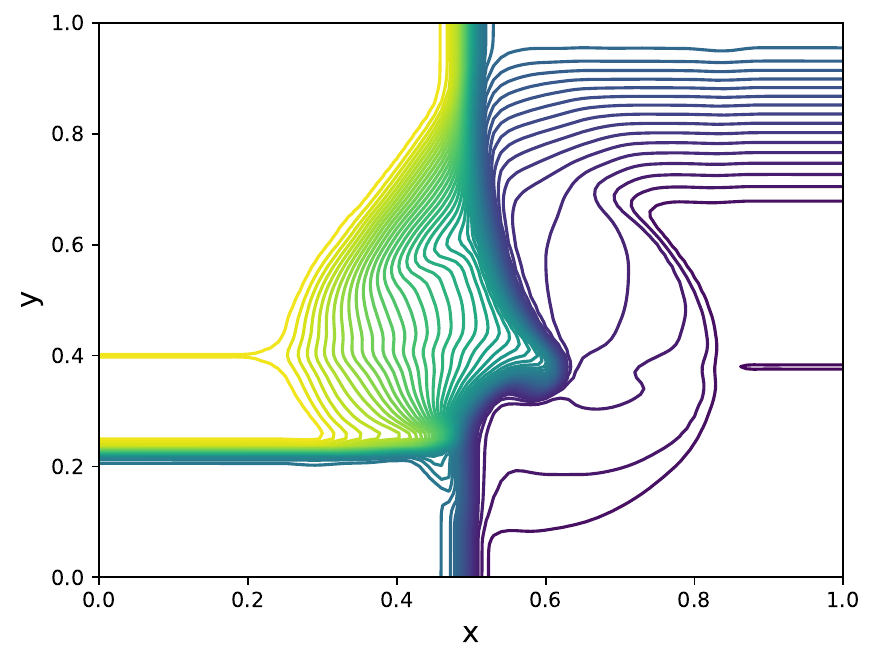} 
\caption{WENO-Z}
\label{fig:R19_WENOZ}
\end{subfigure}
\begin{subfigure}{.32\textwidth}
\centering
\includegraphics[width=\linewidth]{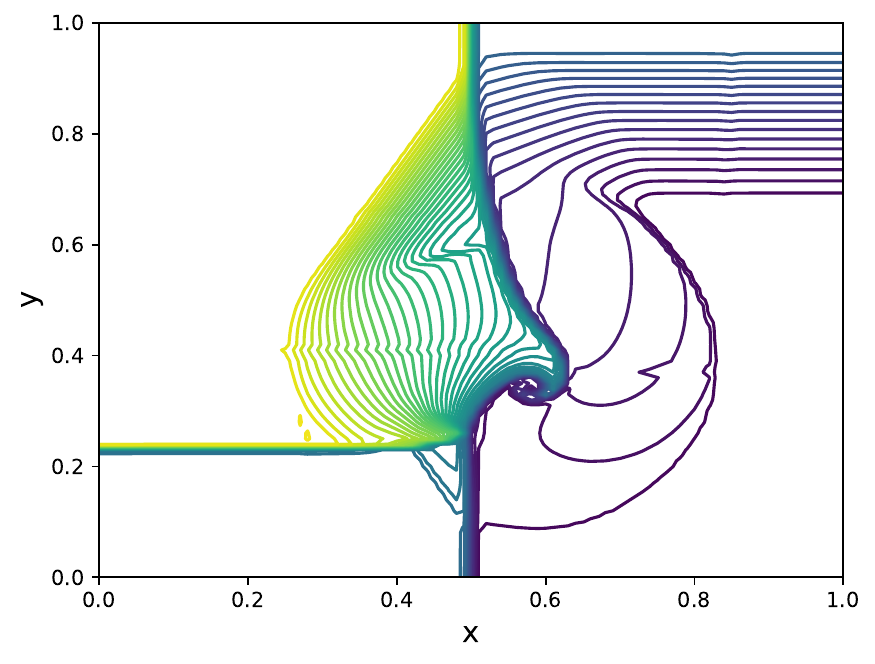} 
\caption{reference solution}
\label{fig:R19_reference}
\end{subfigure}
\begin{subfigure}{.32\textwidth}
\centering
\includegraphics[width=\linewidth]{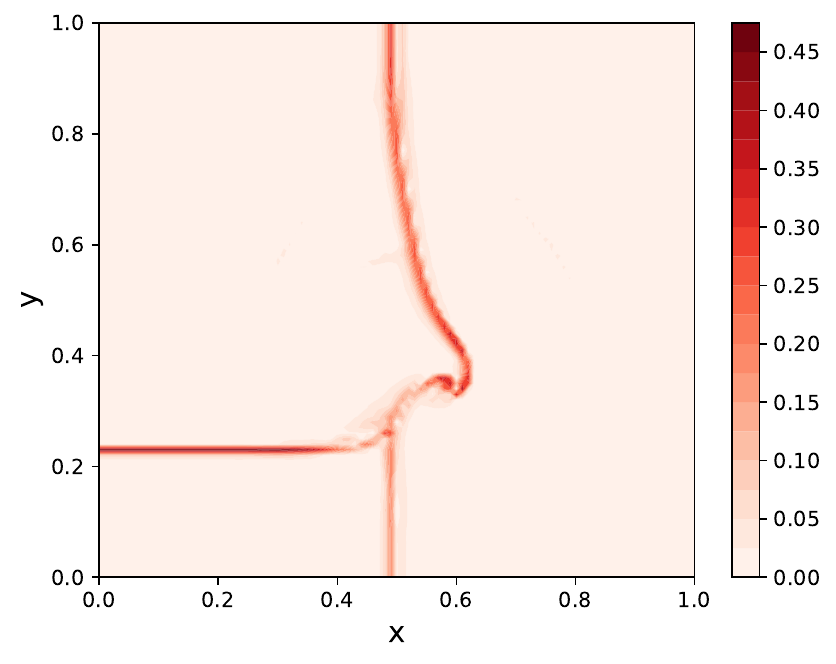} 
\caption{WENO-DS}
\end{subfigure}
\begin{subfigure}{.32\textwidth}
\centering
\includegraphics[width=\linewidth]{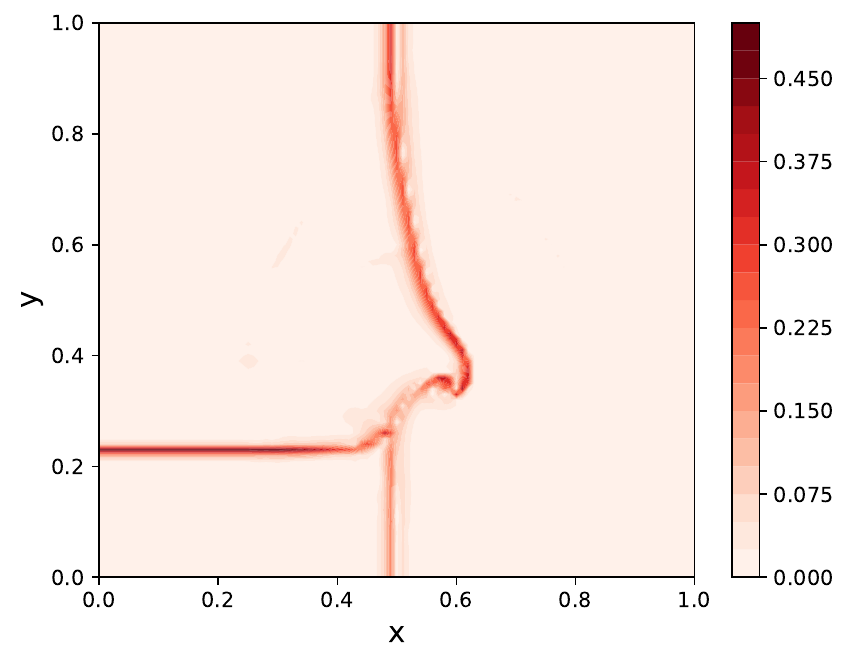} 
\caption{WENO-Z}
\end{subfigure}
\caption{ Density contour plots and absolute pointwise errors for the solution of the Riemann problem with initial Configuration~19, $I \times J = 100 \times 100$, $T=0.3$.} \label{fig:Riemann_19}
\end{figure}

As one can see, we obtain significant error improvements with both methods. 
Comparing both methods, even better results are obtained when the CNN was trained on a data set from Section~\ref{sec:S5.2} on a configuration with four shock waves. 
Compared to the Table~\ref{tab:Riemann_3}, the improvement for Configuration~3 is smaller but still significant. 
However, the method is able to generalize much better to unknown configurations. 
For example, for Configuration~14, we obtain an average improvement rate of $1.30$ for all four variables.
In addition, we use WENO-DS (C3c) to illustrate the density contour plots and absolute pointwise errors in Figures~\ref{fig:Riemann_6}, \ref{fig:Riemann_8}, and \ref{fig:Riemann_19}. 
Here we also show the difference from Configuration~3, with which the model was actually trained.

\begin{figure}[h!] 
\centering
\begin{subfigure}{0.24\textwidth}
    \includegraphics[width=\textwidth]{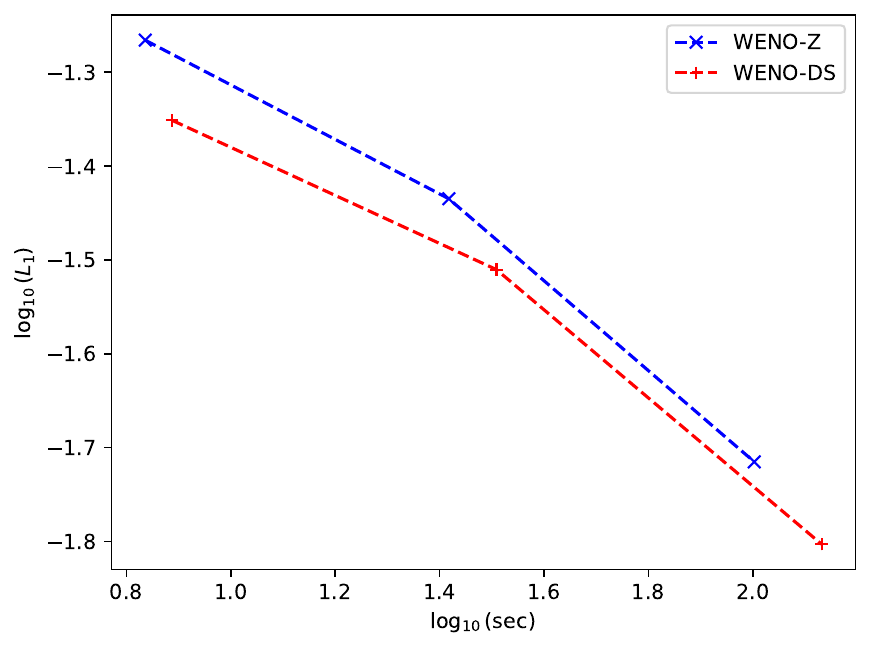}
    \caption{Configuration 3}
    \label{fig:Riemann_3c_cost}
\end{subfigure}
\begin{subfigure}{0.24\textwidth}
     \includegraphics[width=\textwidth]{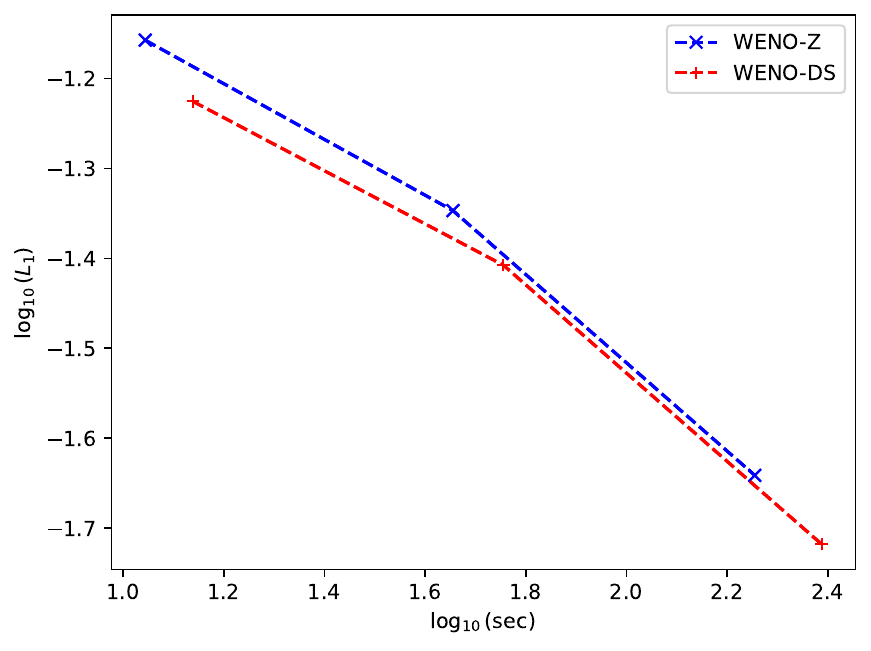}
    \caption{Configuration 6}
        \label{fig:Riemann_6_cost}
\end{subfigure}
\begin{subfigure}{0.24\textwidth}
     \includegraphics[width=\textwidth]{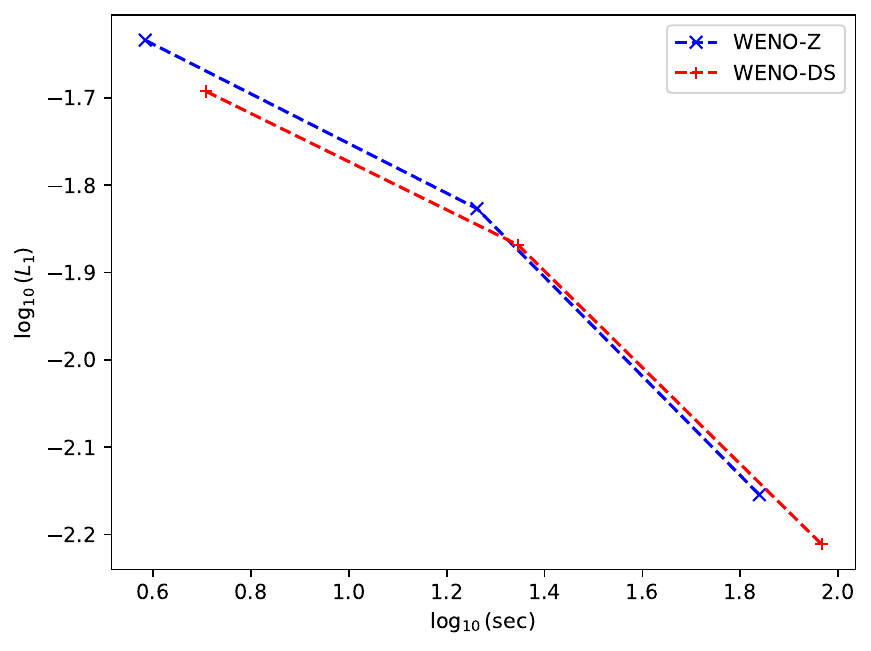}
    \caption{Configuration 8}
        \label{fig:Riemann_8_cost}
\end{subfigure}
\begin{subfigure}{0.24\textwidth}
    \includegraphics[width=\textwidth]{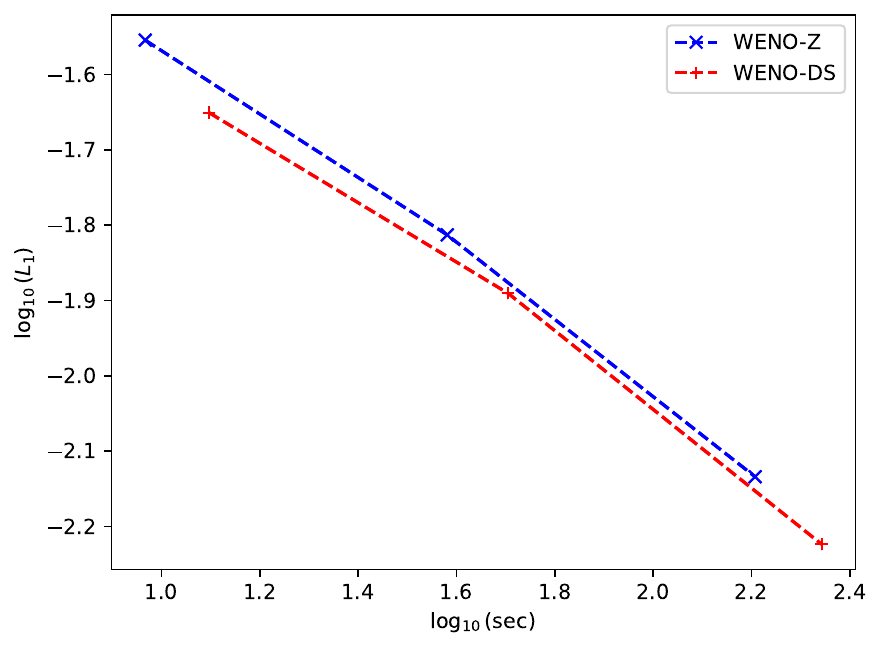}
    \caption{Configuration 11}
    \label{fig:Riemann_11_cost}
\end{subfigure}
\bigskip
\begin{subfigure}{0.24\textwidth}
     \includegraphics[width=\textwidth]{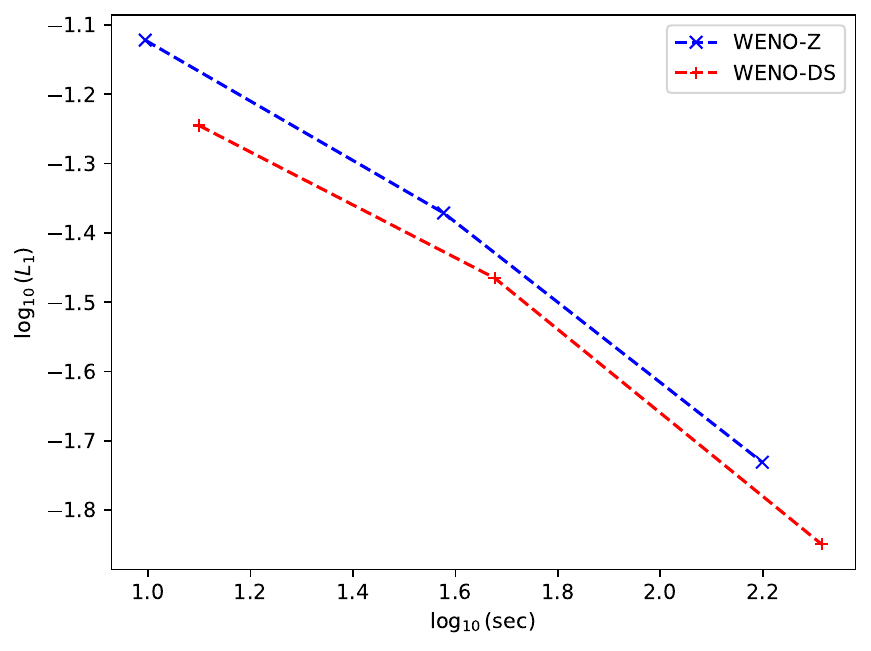}
    \caption{Configuration 14}
        \label{fig:Riemann_14_cost}
\end{subfigure}
\begin{subfigure}{0.24\textwidth}
     \includegraphics[width=\textwidth]{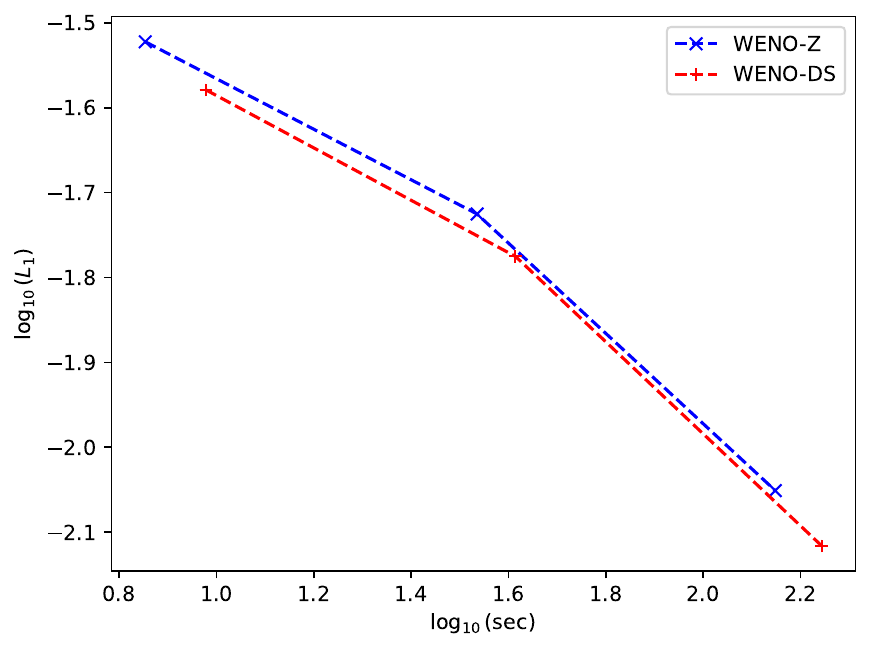}
    \caption{Configuration 18}
        \label{fig:Riemann_18_cost}
\end{subfigure}
\begin{subfigure}{0.24\textwidth}
     \includegraphics[width=\textwidth]{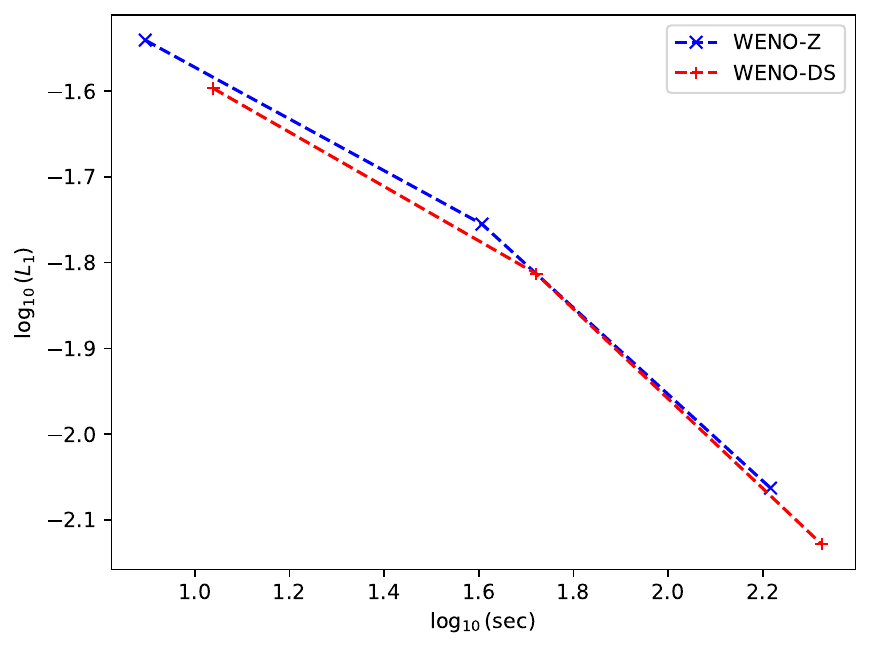}
    \caption{Configuration 19}
        \label{fig:Riemann_19_cost}
\end{subfigure}
\captionsetup{font={footnotesize}}
\caption{Comparison of computational cost against $L_1$-error of the solution of Riemann problem with various initial configurations using WENO-Z and WENO-DS (C3c) methods.}
\label{fig:Riemann_all_costs}
\end{figure}

The WENO-DS (C3c) method achieves large error improvements not only for problems from the same configuration, 
but also for problems from significantly different configurations. 
Since we used a larger CNN, the question is what is the actual numerical cost of these improvements? 
We illustrate the computational costs in Figure~\ref{fig:Riemann_all_costs}. 
As can be seen from the shift of the red dots to the right, the method involves larger computational costs. 
However, it is still more effective or not worse than the original method in most cases.
We would like to emphasize that here we are comparing results with significantly different initial problems 
than those on which the method was actually trained.
Machine learning models are generally not expected to give much better results on unseen problems.

\begin{table}[h] 
 \centering
  \begin{minipage} {.95\linewidth}
   \centering
    \scalebox{0.55}{\begin{tabular}{|c|c|c|c|c|c|c|c|c|c|c|c|c|}
    \hline
    \multicolumn{1}{|c|}{}&\multicolumn{3}{|c|}{\  \textbf{Configuration 3}}&\multicolumn{3}{|c|}{\ Configuration 6} &\multicolumn{3}{|c|}{\ Configuration 8}&\multicolumn{3}{|c|}{\ Configuration 11}\\
    \hline
        & \textbf{WENO-Z} &   \textbf{WENO-DS} &  \textbf{ratio} & WENO-Z &  WENO-DS & ratio& WENO-Z &  WENO-DS & ratio& WENO-Z &  WENO-DS & ratio\\
    \hline
    \hline
 $\rho$  &   \textbf{0.019232} &   \textbf{0.015033} &   \textbf{1.28}  &  0.038616 &  0.032696 &  1.18  &  0.005711 &  0.004975 &  1.15 &  0.007792 &  0.006316 &  1.23  \\  \hline
$u$    &  \textbf{ 0.019588} &   \textbf{0.016359} &   \textbf{1.20}  &  0.019662 &  0.016144 &  1.22  &  0.008488 &  0.007396 &  1.15 &  0.008003 &  0.006487 &  1.23  \\  \hline
$v$   &   \textbf{0.019588} &   \textbf{0.016359} &   \textbf{1.20}  &  0.022582 &  0.018951 &  1.19  &  0.008488 &  0.007396 &  1.15 &  0.007692 &  0.006282 &  1.22 \\     \hline
 $p$   &   \textbf{0.018666} &   \textbf{0.015214} &   \textbf{1.23}  &  0.010525 &  0.008821 &  1.19  &  0.005350 &  0.004844 &  1.10  &  0.005883 &  0.004813 &  1.22 \\  \hline
\end{tabular} }

 \vspace{10pt}
\scalebox{0.55}{
\begin{tabular}{|c|c|c|c|c|c|c|c|c|c|}
    \hline
    \multicolumn{1}{|c|}{}&\multicolumn{3}{|c|}{\ Configuration 14}&\multicolumn{3}{|c|}{\ Configuration 18} &\multicolumn{3}{|c|}{\ Configuration 19}\\
    \hline
        & WENO-Z &  WENO-DS & ratio & WENO-Z &  WENO-DS & ratio& WENO-Z &  WENO-DS & ratio\\
    \hline
    \hline
 $\rho$    &  0.013169 &  0.010333 &  1.27 &  0.014918 &  0.012519 &  1.19  &  0.014844 &  0.012390 &  1.20 \\    \hline
$u$     &  0.004835 &  0.003732 &  1.30  &  0.003534 &  0.003063 &  1.15  &  0.003749 &  0.003339 &  1.12 \\     \hline
$v$       &  0.021299 &  0.016512 &  1.29  &  0.010315 &  0.009077 &  1.14  &  0.009891 &  0.008641 &  1.14 \\    \hline
 $p$       &  0.034996 &  0.026008 &  1.35  &  0.006795 &  0.005961 &  1.14 &  0.006123 &  0.005393 &  1.14 \\   \hline
\end{tabular} }
  \end{minipage}
    \caption{Comparison of $L_1$ error of WENO-Z and WENO-DS (C3c) methods for the solution of the Euler system with various initial configurations, $I \times J = 100 \times 100$.}
\label{tab:Riemann_C3c}
\end{table}

\begin{table}[h] 
 \centering
  \begin{minipage} {.95\linewidth}
   \centering
    \scalebox{0.55}{\begin{tabular}{|c|c|c|c|c|c|c|c|c|c|c|c|c|}
    \hline
    \multicolumn{1}{|c|}{}&\multicolumn{3}{|c|}{\  \textbf{ \textbf{Configuration 16}}}&\multicolumn{3}{|c|}{\ Configuration 6} &\multicolumn{3}{|c|}{\ Configuration 8}&\multicolumn{3}{|c|}{\ Configuration 11}\\
    \hline
        & \textbf{WENO-Z} &   \textbf{WENO-DS} &  \textbf{ratio} & WENO-Z &  WENO-DS & ratio& WENO-Z &  WENO-DS & ratio& WENO-Z &  WENO-DS & ratio\\
    \hline
    \hline
 $\rho$  &  \textbf{ 0.004834} &   \textbf{0.004127} &   \textbf{1.17}  &  0.038616 &  0.036329 &  1.06  &  0.005711 &  0.004777 &  1.20 &  0.007792 &  0.007695 &  1.01 \\      \hline
$u$    &   \textbf{0.005989} &   \textbf{0.004981} &   \textbf{1.20}  &  0.019662 &  0.019575 &  1.00  &  0.008488 &  0.007056 &  1.20  &  0.008003 &  0.007824 &  1.02 \\      \hline
$v$   &   \textbf{0.006609} &   \textbf{0.005776} &   \textbf{1.14}  &  0.022582 &  0.019974 &  1.13  &  0.008488 &  0.007056 &  1.20  &  0.007692 &  0.007482 &  1.03 \\        \hline
 $p$   &   \textbf{0.004236} &   \textbf{0.003663} &   \textbf{1.16}  &  0.010525 &  0.010216 &  1.03  &  0.005350 &  0.004624 &  1.16  &  0.005883 &  0.006295 &  0.93 \\     \hline
\end{tabular} }

 \vspace{10pt}
\scalebox{0.55}{
\begin{tabular}{|c|c|c|c|c|c|c|c|c|c|}
    \hline
    \multicolumn{1}{|c|}{}&\multicolumn{3}{|c|}{\ Configuration 14}&\multicolumn{3}{|c|}{\ Configuration 18} &\multicolumn{3}{|c|}{\ Configuration 19}\\
    \hline
       & WENO-Z &  WENO-DS & ratio & WENO-Z &  WENO-DS & ratio& WENO-Z &  WENO-DS & ratio\\
    \hline
    \hline
 $\rho$     &  0.013169 &  0.011718 &  1.12  &  0.014918 &  0.013447 &  1.11  &  0.014844 &  0.013198 &  1.12 \\    \hline
$u$      &  0.004835 &  0.004042 &  1.20  &  0.003534 &  0.002975 &  1.19  &  0.003749 &  0.003256 &  1.15 \\   \hline
$v$      &  0.021299 &  0.020330 &  1.05  &  0.010315 &  0.009302 &  1.11  &  0.009891 &  0.008796 &  1.12 \\  \hline
 $p$        &  0.034996 &  0.036038 &  0.97  &  0.006795 &  0.006535 &  1.04  &  0.006123 &  0.005752 &  1.06 \\ \hline
\end{tabular} }
  \end{minipage}
    \caption{Comparison of $L_1$ error of WENO-Z and WENO-DS (C16c) methods for the solution of the Euler system with various initial configurations, $I \times J = 100 \times 100$.}
\label{tab:Riemann_C16c}
\end{table}

\section{Conclusion} \label{sec:S6}
In this paper, we introduced a novel approach, WENO-DS, which leverages the power of deep learning to enhance the performance of the well-established Weighted Essentially Non-Oscillatory (WENO) scheme in the context of solving hyperbolic conservation laws, particularly exemplified by the two-dimensional Euler equations of gas dynamics. 
By seamlessly integrating deep learning techniques into the WENO algorithm,
we have successfully improved the accuracy of numerical solutions, particularly in regions near abrupt shocks. 
Unlike previous attempts at incorporating deep learning into numerical methods, this approach stands out by eliminating the need for additional post-processing steps, ensuring consistency throughout.

This study demonstrates the superiority of the WENO-DS approach \linebreak through an extensive examination of various test problems, including scenarios featuring shocks and rarefaction waves.
The results consistently showcase the newfound capabilities of the approach, outperforming traditional fifth-order WENO schemes, especially when dealing with challenges like excessive diffusion or overshooting around shocks.

The introduction of machine learning into the realm of solving partial differential equations (PDEs) has brought about promising improvements in numerical methods. 
However, it is crucial to strike a balance between these data-driven insights and the foundational mathematical principles underpinning the numerical scheme. 
This study successfully maintains this equilibrium, building upon the physical principles of the Euler equations while incorporating deep learning enhancements.

In summary, the WENO-DS approach represents a significant advancement in the field of numerical methods for hyperbolic conservation laws, where the incorporation of deep learning techniques has not only enhanced the accuracy but also improved the qualitative behavior of solutions, both in smooth regions and near discontinuities. 
This research paves the way for future developments in the intersection of traditional numerical methods and machine learning, offering a promising direction for further advancements in solving complex PDEs like the Euler equations.



 \bibliographystyle{elsarticle-num} 
 \bibliography{cas-refs}





\end{document}